\pdfoutput=1
\documentclass[11pt]{article}
\usepackage{bbm}
\usepackage{amsmath}
\usepackage{amsthm}
\usepackage{amssymb} 
\usepackage{mathrsfs}
\usepackage[dvipsnames]{xcolor}
\usepackage{graphicx}
\usepackage{authblk}
\usepackage[backref=none,hypertexnames=false,colorlinks=true,urlcolor=blue,linkcolor=blue,citecolor=blue]{hyperref}
\usepackage[numbers,comma,square,sort&compress]{natbib}
\usepackage[letterpaper,text={6.5in,9in},centering]{geometry}
\usepackage{tabularx}
\usepackage{booktabs} 
\usepackage{algorithm}
\usepackage{array} 
\usepackage{subcaption} 
\providecommand{\keywords}[1]{\textbf{Keywords:} #1}

\renewcommand{\baselinestretch}{1.1}

\graphicspath{{eps/}{pdf/}}

\makeatletter\@addtoreset{equation}{section}\makeatother

\DeclareMathOperator*{\argmin}{arg\,min}

\begin{document}

\title{Weighted Birkhoff Averages Accelerate Data-Driven Methods}
\author[1]{Maria Bou-Sakr-El-Tayar}
\author[2]{Jason J. Bramburger}
\author[3]{Matthew J. Colbrook}

\affil[1]{\small Department of Computer Science and Software Engineering, Concordia University, Montr\'eal, QC, Canada}

\affil[2]{\small Department of Mathematics and Statistics, Concordia University, Montr\'eal, QC, Canada}

\affil[3]{\small Department of Applied Mathematics and Theoretical Physics, Cambridge University, Cambridge, United Kingdom}

\date{}
\maketitle

\begin{abstract}
Many data-driven algorithms in dynamical systems rely on ergodic averages that converge painfully slowly. One simple idea changes this: taper the ends. Weighted Birkhoff averages can converge much faster (sometimes superpolynomially, even exponentially) and can be incorporated seamlessly into existing methods. We demonstrate this with five weighted algorithms: weighted Dynamic Mode Decomposition (wtDMD), weighted Extended DMD (wtEDMD), weighted Sparse Identification of Nonlinear Dynamics (wtSINDy), weighted spectral measure estimation, and weighted diffusion forecasting. Across examples ranging from fluid flows to El Ni\~no data, the message is clear: weighting costs nothing, is easy to implement, and often delivers markedly better results from the same data.
\end{abstract}


\section{Introduction}

Many algorithms in dynamical systems rely on ergodic averages.
Given data from a single trajectory $x_{n+1} = f(x_n)$, one forms sums such as the Birkhoff average
\begin{equation}\label{BirkhoffAverage}
    B_N(g)(x)=\frac{1}{N}\sum_{n = 0}^{N-1} g(f^n(x))   
\end{equation}
hoping that as $N$ grows, the average approaches an expectation.
This idea lies at the heart of data sampling assumptions for diffusion forecasting \cite{giannakis2019data,berry2015nonparametric,thiede2019galerkin,giannakis2020extraction,berry2015nonparametric2}, Koopman operator analysis \cite{korda2018convergence,bramburger2024auxiliary}, Dynamic Mode Decomposition \cite{korda2018convergence,colbrook2023beyond,bramburger2024auxiliary},
and other approaches such as system identification \cite{russo2024convergence} and power spectra \cite{colbrook2024rigorous,colbrook2025rigged}.

However, the classical Birkhoff average converges painfully slowly. In typical cases, the error decays only like $1/N$ \cite{kachurovskii1996rate}, and for certain $g$, it can be arbitrarily slower \cite{krengel1978speed}. Long trajectories are therefore needed to obtain reliable results. These slow rates set the pace for data-driven methods built upon them, dictating how much data is required for accurate approximation.

The remedy is simple but powerful: replace the uniform average with a weighted one \cite{tong2025weighted,duignan2023distinguishing,das2018super,das2017quantitative,tong2024exponential}. By tapering the ends of the sum, edge effects are suppressed and convergence can improve dramatically, sometimes to super-polynomial or even exponential rates. Weighted ergodic averages thus provide sharper estimates from the same data and improve a wide range of spectral and operator-learning methods.

In what follows, we review the basic idea of weighted Birkhoff averages and illustrate, through five examples, how they enhance the accuracy of data-driven computations.

\section{Weighted Birkhoff averages}\label{sec:WBAs}

Consider a measure space $(\mathcal{X},\mathcal{A},\mu)$ and a measure-preserving map $f:\mathcal{X}\to \mathcal{X}$. The measure $\mu$ is said to be {\em ergodic} if every set $S \in \mathcal{A}$ satisfying $f^{-1}(S)= S$ has measure $\mu(S) \in \{0,1\}$. Birkhoff’s ergodic theorem \cite[Theorem~4.5.5]{brin2002introduction} then states that if $\mu$ is ergodic and $g \in L^1(\mathcal{X},\mu)$, the time average  $B_N(g)(x)$ presented in \eqref{BirkhoffAverage} converges for $\mu$-almost every $x \in \mathcal{X}$ 
to the space average
$$
    \langle g\rangle := \int_\mathcal{X} g\ \mathrm{d}\mu
$$
as $N \to \infty$. To improve this convergence, we turn to the theory of weighted Birkhoff averages.   

\paragraph{Weighted averages.} 
Define the class of smooth weight functions
\begin{equation}\label{WeightClass}
    \mathcal{W} = \bigg\{w \in C^\infty([0,1],\mathbb{R}_+):\ \int_0^1 w(x)\mathrm{d}x = 1,\ \frac{\mathrm{d}^{i}w}{\mathrm{d}x^{i}}(0) = \frac{\mathrm{d}^{i}w}{\mathrm{d}x^{i}}(1) = 0,\ \forall i = 0,1,2,\dots\bigg\}.
\end{equation}
Each $w\in\mathcal{W}$ is a smooth ``bump'' function whose derivatives vanish at the boundaries.
A standard example is
\begin{equation}\label{BumpFunction}
    w(x) = \begin{cases}
        C\mathrm{e}^{\frac{-1}{x(1-x)}}, & x \in (0,1), \\
        0, & x = 0,1,
    \end{cases}
\end{equation}
where the constant $C$ ensures that $\int_0^1 w(x)\mathrm{d}x = 1$. For such a weight, we define the \textit{weighted Birkhoff average} as
\begin{equation}\label{WeightedBirkhoffAverage}
    WB_N(g)(x) := \frac{1}{\alpha_N} \sum_{n = 0}^{N-1} w(n/N)g(f^n(x)), \quad \mathrm{where} \quad \alpha_N := \sum_{n = 0}^{N-1} w(n/N).
\end{equation}
The weight $w$ damps the ends of the sum, giving greater emphasis to terms near the centre $n \sim N/2$, while the uniform average $B_N$ treats all terms equally. If $\mu$ is ergodic, then for any $w \in \mathcal{W}$ and $g \in L^1(\mathcal{X},\mu)$ we have \cite[Proposition~2]{duignan2023distinguishing}
\begin{equation}\label{WeightedConvergence}
    \lim_{N \to \infty} WB_N(g)(x) = \langle g\rangle 
\end{equation}
for $\mu$-almost every $x \in X$. Thus, both $WB_N(g)(x)$ and $B_N(g)(x)$ converge to the same limit, but possibly at very different speeds.

As summarised in \cite{tong2025weighted}, the convergence rate can be greatly improved beyond the unweighted $\mathcal{O}(1/N)$ rate. If $f$ is periodic, the rate of convergence in \eqref{WeightedConvergence} is exponential. If $f$ is quasiperiodic and $g \in C^\infty$, the rate is $\mathcal{O}(N^{-m})$ for any $ m > 0$; if $g$ is analytic, it again becomes exponential.
For chaotic dynamics, no general analytic rate is known, but empirical evidence (including that below) suggests that the weighted and unweighted averages share the same $\mathcal{O}(1/N)$ behaviour.

{\color{black}For smooth quasiperiodic dynamics, the slow convergence of Birkhoff averages is mainly a boundary artifact caused by the sharp truncation of the time sum. The average behaves like a low-order quadrature applied to a smooth periodic function, with errors dominated by endpoint effects rather than interior smoothness. Weighted averages replace this abrupt cutoff with a smooth window whose derivatives vanish at the endpoints, suppressing oscillatory Fourier modes by repeated integration by parts. Each non-constant frequency is then damped faster than any algebraic power of the averaging length, while the zero mode (the desired average) is unchanged. In effect, the weighting upgrades the time average to a high-order quadrature rule, yielding super-polynomial convergence for smooth quasiperiodic systems with Diophantine frequencies. This mechanism relies on the regular spectral structure of quasiperiodicity and does not generally extend to chaotic dynamics, where broadband correlations dominate.}

\begin{figure}[t]
    \center
    \begin{minipage}[t]{0.32\textwidth}
        \centering
        \includegraphics[width = \textwidth]{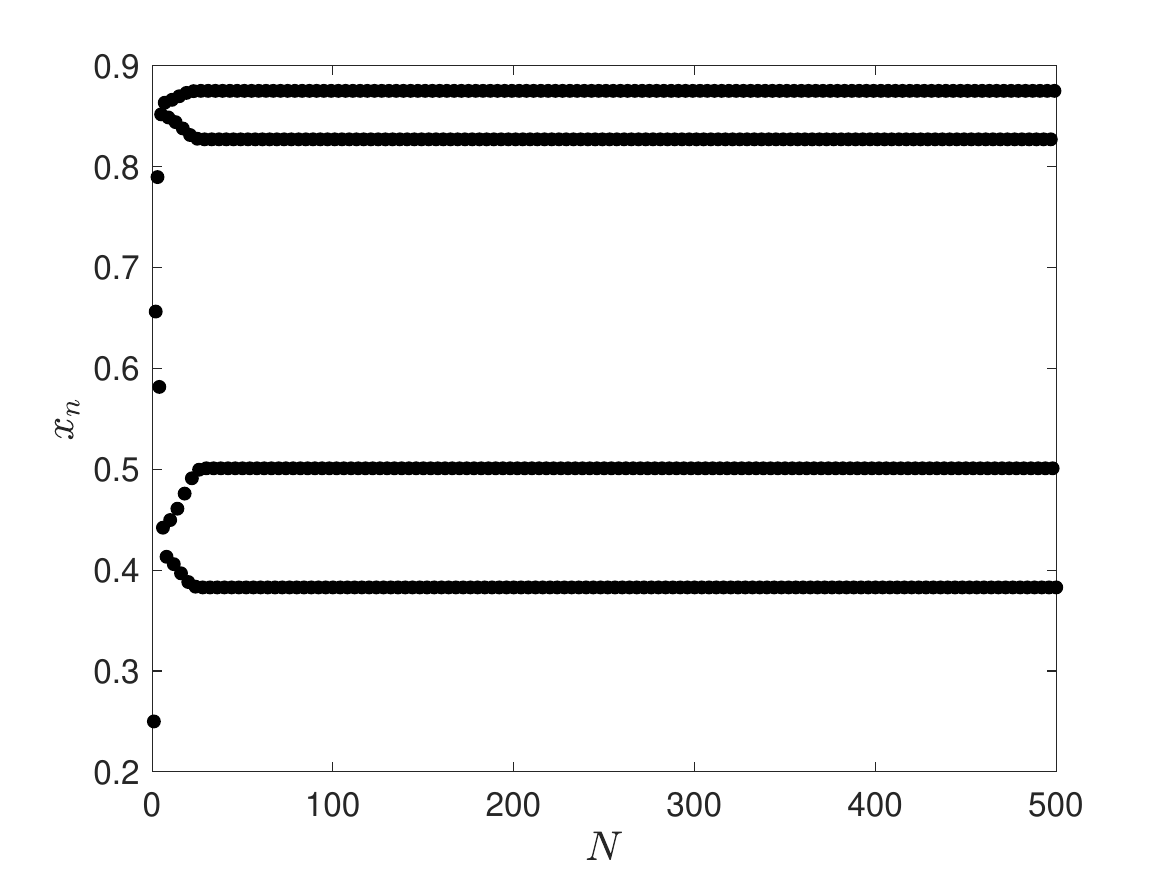}
    \end{minipage}
     \begin{minipage}[t]{0.32\textwidth}
        \centering
        \includegraphics[width = \textwidth]{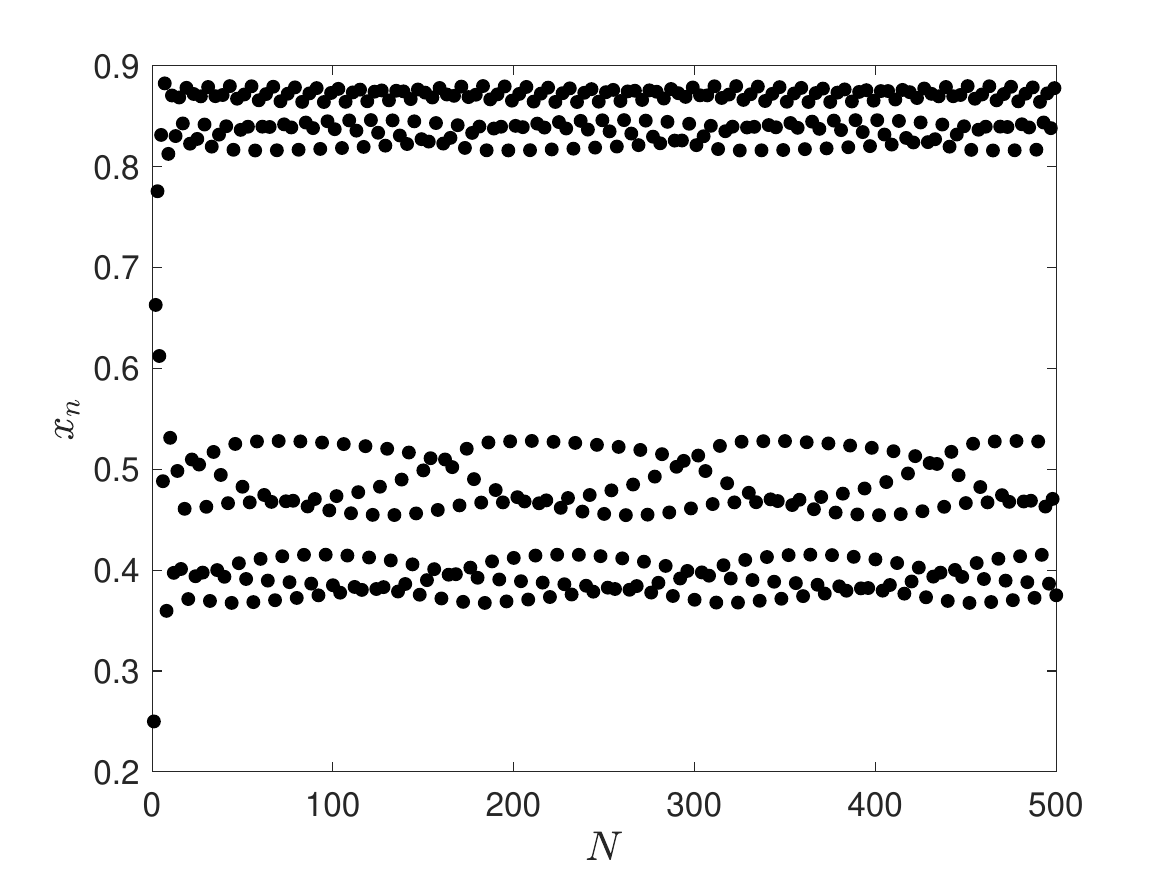}
    \end{minipage}
    \begin{minipage}[t]{0.32\textwidth}
        \centering
        \includegraphics[width = \textwidth]{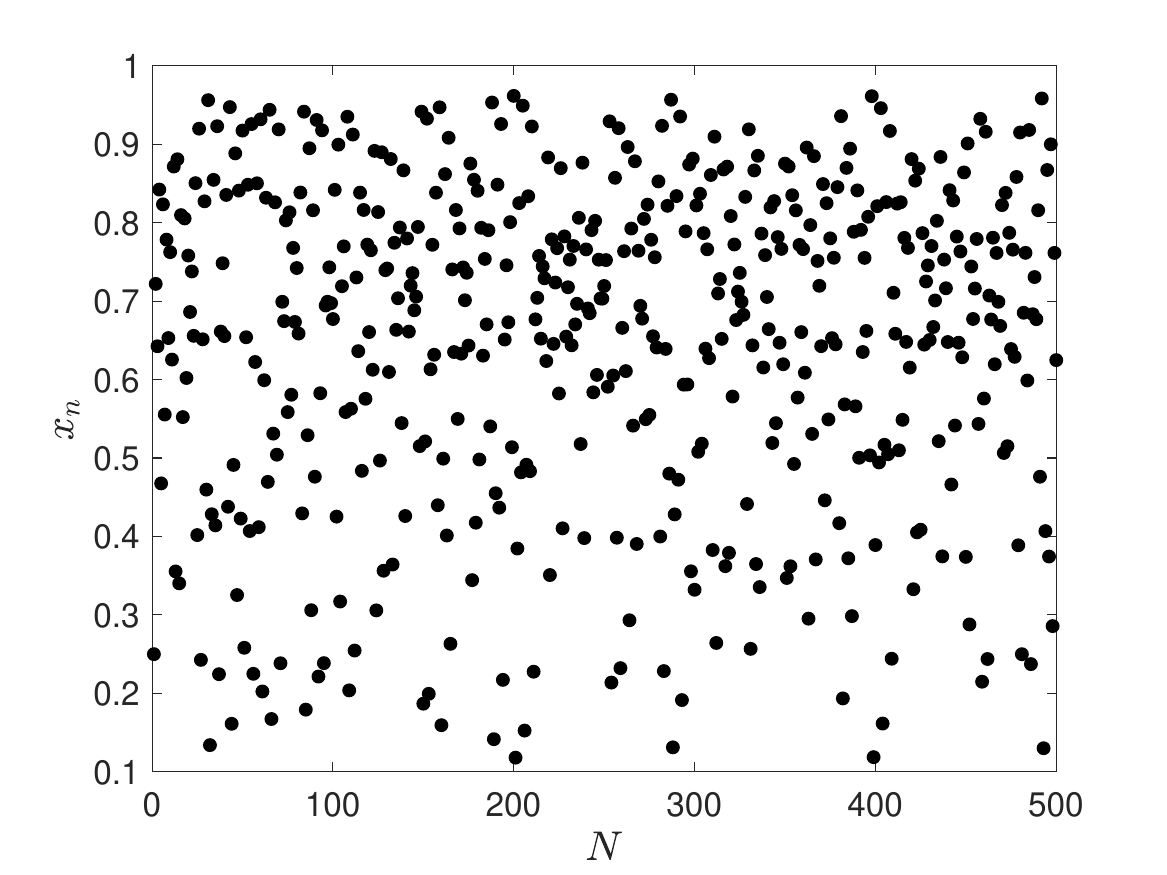}
    \end{minipage}
    \begin{minipage}[t]{0.32\textwidth}
        \centering
        \includegraphics[width = \textwidth]{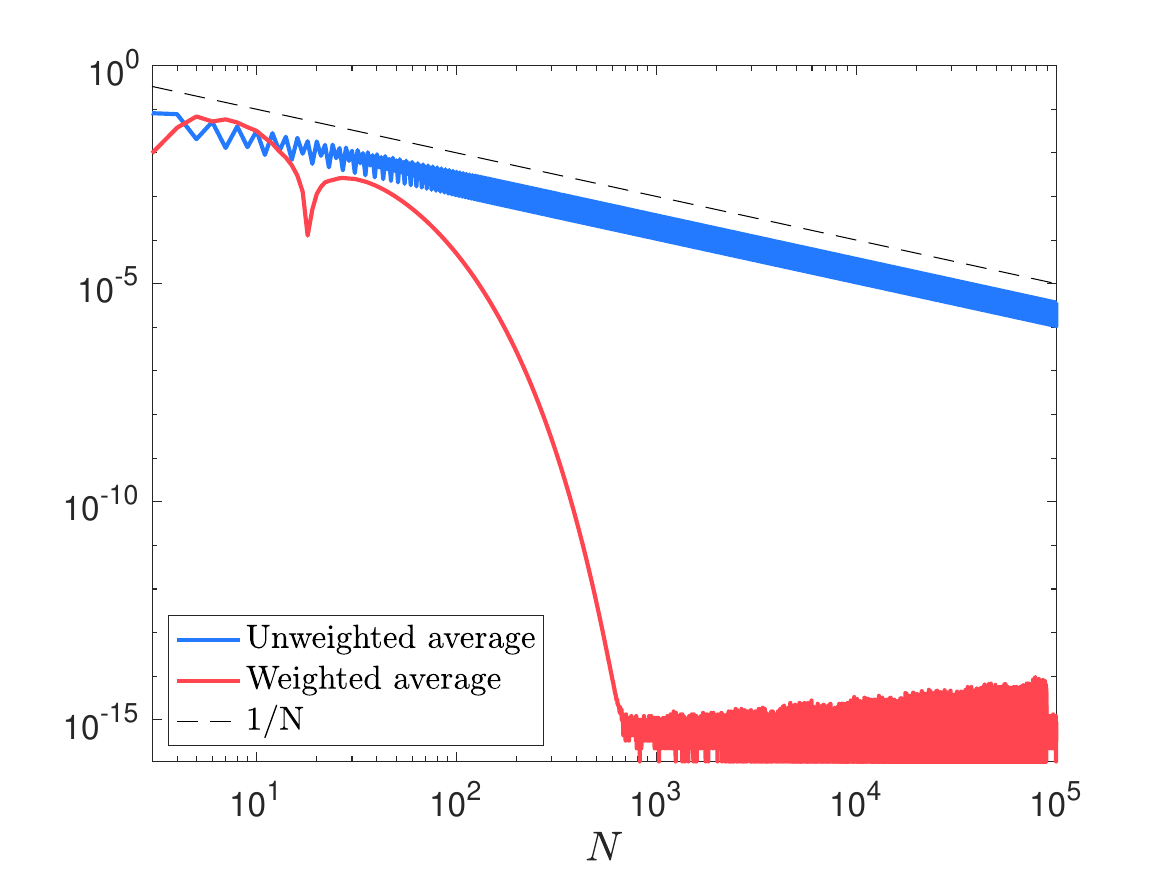}
    \end{minipage}
    \begin{minipage}[t]{0.32\textwidth}
        \centering
        \includegraphics[width = \textwidth]{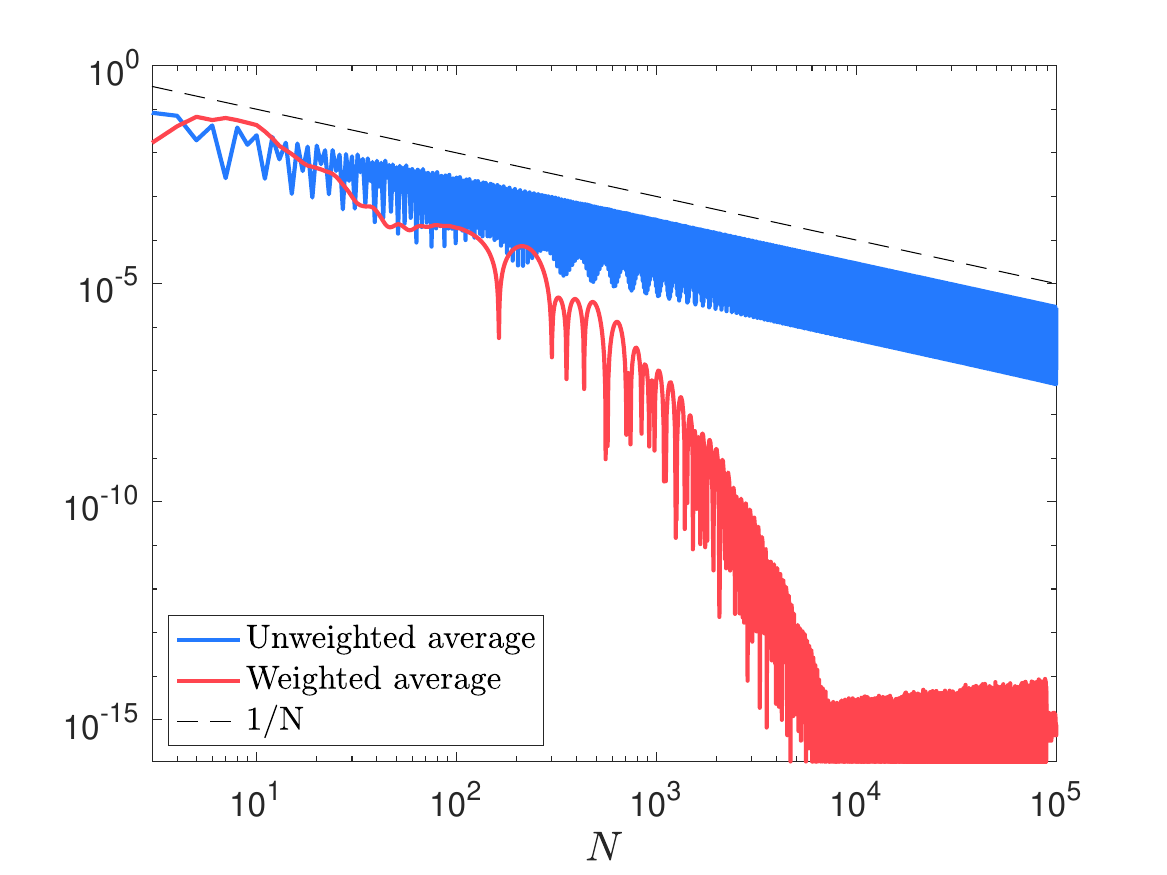}
    \end{minipage}
    \begin{minipage}[t]{0.32\textwidth}
        \centering
        \includegraphics[width = \textwidth]{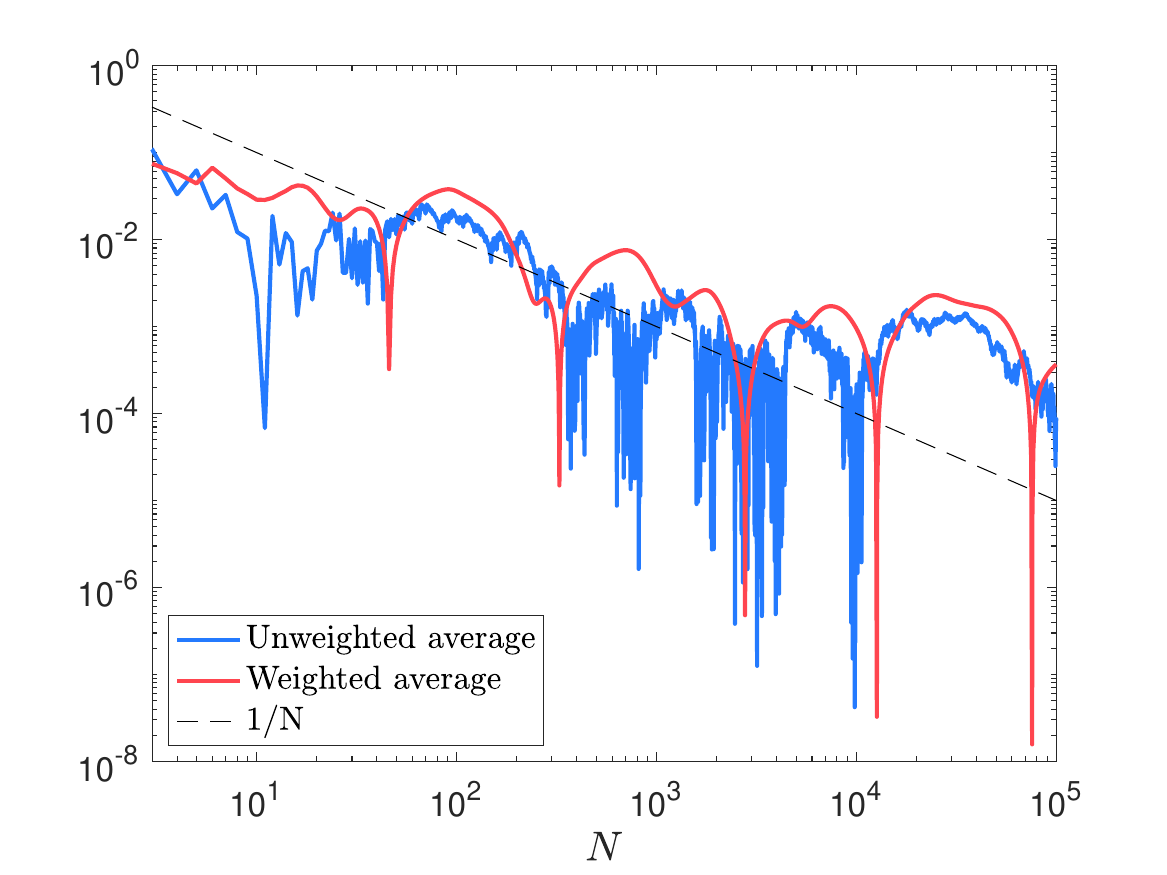}
    \end{minipage}
    \caption{Dynamics (top) and averages (bottom) of the driven logistic map \eqref{DrivenLogistic}. Weighted (red) and unweighted (blue) time averages of the state variable $x$ are presented as a log-log plot of absolute errors against the average with $N = 10^8$. Parameter values are $\varepsilon = 0$ (left, periodic), $\varepsilon = 0.01$ (middle, quasiperiodic), and $\varepsilon = 0.1$ (right, chaotic).}
    \label{fig:Logistic}
\end{figure}

\paragraph{Demonstration on the driven logistic map.} To illustrate, consider the driven logistic map \cite{negi2000bifurcations}:
\begin{equation}\label{DrivenLogistic}
        x_{n+1} = 3.5(1 + \varepsilon \cos(2\pi \theta_n))x_n(1-x_n) \qquad
        \theta_{n+1} = \theta_n + \sqrt{2} \mod 1,
\end{equation}
where the parameter $\varepsilon \geq 0$ controls the strength of a quasiperiodic forcing.
We examine time averages of $g(x,\theta) = x$ for $\varepsilon = 0$ (periodic), $\varepsilon = 0.01$ (quasiperiodic) and $\varepsilon = 0.1$ (chaotic).  

The top row of Figure \ref{fig:Logistic} shows the iterates $x_n$ for these three regimes, all with $x_0 = 0.25$. The bottom row plots, on a log-log scale, the absolute error between the running average and the value at $N = 10^8$. The unweighted average converges at roughly the $\mathcal{O}(1/N)$ rate for all cases.
For non-chaotic dynamics ($\varepsilon = 0,0.01$), the weighted averages converge exponentially to machine precision; the small oscillations between $10^{-14}$ and $10^{-16}$ reflect floating-point noise.
When $\varepsilon = 0.1$ the map is chaotic, and the weighted and unweighted averages converge at nearly the same rate. Thus, weighting never slows convergence and can dramatically accelerate it whenever the underlying dynamics are not chaotic.

\section{Five Examples of Weighted Methods}

We build on several well-established data-driven methods for analysing dynamical datasets.
Background and fundamentals can be found in \cite{bramburger2024databook,colbrook_book,brunton2022data}.
Here, we show how these methods can be strengthened by incorporating weighted Birkhoff averages.
When a specific weight $w \in \mathcal{W}$ is needed for illustration, we use the bump function defined in~\eqref{BumpFunction}. All examples that follow can be reproduced with the freely available code at
 \href{https://github.com/jbramburger/weighted_methods}{https://github.com/jbramburger/weighted\_methods}.

\subsection{Example Method I: Dynamic mode decomposition}

Dynamic Mode Decomposition (DMD) is a data-driven technique for extracting structure from sequential data \cite{tu2014dynamic,schmid2010dynamic}; see also the review~\cite{colbrook2024multiverse}.
Given snapshots $X_1,X_2,\dots,X_N,X_{N+1} \in \mathbb{R}^d$, for some $N,d \geq 1$, the goal is to identify a best-fit linear map for $X_n \mapsto X_{n+1}$. In its simplest form, DMD seeks a matrix $A \in \mathbb{R}^{d\times d}$ so that $X_{n+1} \approx AX_n$ for all $n = 1,\dots, N$ by minimizing $\|\mathbb{Y} - A\mathbb{X}\|_F$, where $\|\cdot\|_F$ is the Frobenius norm and
\begin{equation}\label{SnapshotMatrices}
        \mathbb{X} = \begin{bmatrix}
            X_1 & \cdots & X_N
        \end{bmatrix},\qquad  \mathbb{Y} = \begin{bmatrix}
            X_2 & \cdots & X_{N+1}
        \end{bmatrix} \in\mathbb{R}^{d \times N}.
\end{equation}
The optimal solution is the DMD matrix $A = \mathbb{Y}\mathbb{X}^\dagger$, where $\dagger$ denotes the Moore--Penrose pseudoinverse. {\color{black} The Moore--Penrose pseudoinverse can be computed from the singular value decomposition of a matrix $\mathbb{X} = U\Sigma V^\top$ as $\mathbb{X}^\dagger = V\Sigma^\dagger U^\top$ with the pseudoinverse of the rectangular diagonal matrix $\Sigma$, denoted $\Sigma^\dagger$, obtained by taking the reciprocal of each non-zero element of $\Sigma$ and leaving the zeros in place.} The eigenvectors of $A$ (DMD modes) represent coherent spatial structures that grow, decay, or oscillate linearly in time according to their associated eigenvalues.

{\color{black}Using the pseudoinverse identity $\mathbb{X}^\dagger = \mathbb{X}^\top(\mathbb{X}\mathbb{X}^\top)^\dagger$}, this matrix can be written equivalently as
$$
    A = \mathbb{Y}\mathbb{X}^\dagger = \bigg(\frac{1}{N}\mathbb{Y}\mathbb{X}^\top\bigg)\bigg(\frac{1}{N}\mathbb{X}\mathbb{X}^\top\bigg)^\dagger.
$$
which involves the Birkhoff averages
\begin{equation}\label{DMDaverages}
        \frac{1}{N}\mathbb{Y}\mathbb{X}^\top = \frac{1}{N} \sum_{n = 1}^N X_{n+1}X_{n}^\top, \qquad
        \frac{1}{N}\mathbb{X}\mathbb{X}^\top = \frac{1}{N} \sum_{n = 1}^N X_{n}X_{n}^\top.
\end{equation}
Convergence results for DMD as $N \to \infty$ \cite{korda2018convergence,colbrook2023beyond,bramburger2024auxiliary} rely on the convergence of these averages; the limiting matrix represents a projection of the Koopman operator onto the snapshot coordinates (see the following subsection).

\paragraph{Weighted DMD.}
To accelerate convergence, we introduce weighted averages, resulting in the weighted DMD (wtDMD) method, which is particularly useful when data are scarce. Let $w \in \mathcal{W}$ be a weight function from~\eqref{WeightClass},
and define the temporally weighted snapshot matrices
$$
        \mathbb{X}_w := \mathbb{X}W^{1/2}, \qquad \mathbb{Y}_w := \mathbb{Y}W^{1/2}, 
$$
where $W := \mathrm{diag}(w(0),w(1/N),\dots,w((N-1)/N))) \in \mathbb{R}^{N\times N}$. The weighted DMD matrix is
$$
    A_w = \mathbb{Y}_w\mathbb{X}_w^\dagger = \bigg(\frac{1}{\alpha_N}\mathbb{Y}_w\mathbb{X}_w^\top\bigg)\bigg(\frac{1}{\alpha_N}\mathbb{X}_w\mathbb{X}_w^\top\bigg)^\dagger,
$$
where $\alpha_N$ are defined in \eqref{WeightedBirkhoffAverage}. The weighted averages are
\begin{equation}\label{wtDMDaverages}
        \frac{1}{\alpha_N}\mathbb{Y}_w\mathbb{X}_w^\top = \frac{1}{\alpha_N} \sum_{n = 1}^N X_{n+1}X_{n}^\top w((n-1)/N), \qquad
        \frac{1}{\alpha_N}\mathbb{X}_w\mathbb{X}_w^\top = \frac{1}{\alpha_N} \sum_{n = 1}^N X_{n}X_{n}^\top w((n-1)/N).
\end{equation}
Since these converge to the same limit as the unweighted averages in~\eqref{DMDaverages},
the weighted formulation offers potentially faster convergence of the DMD matrix.
Algorithm~\ref{alg:wtDMD} summarises the procedure.

\begin{algorithm}[t]
    \caption{Weighted Dynamic Mode Decomposition (wtDMD)}\label{alg:wtDMD}
    {\bf Input:} Sequential snapshots $\{X_n\}_{n = 1}^{N+1}$ and bump function $w \in \mathcal{W}$.
    \begin{enumerate}
        \item Arrange data into matrices $\mathbb{X}$ and $\mathbb{Y}$ defined in \eqref{SnapshotMatrices} and $W = \mathrm{diag}(w(0),w(1/N),\dots,w((N-1)/N))$. 
        \item Compute the weighted DMD matrix $A_w = (\mathbb{X} W^{1/2})(\mathbb{Y} W^{1/2})^\dagger$.
        \item Perform eigendecomposition $A_wV = \Lambda V$ to identify DMD modes.
    \end{enumerate}
    {\bf Output:} Weighted DMD matrix $A_w$, with eigenvectors $V$ and eigenvalues $\Lambda$.
\end{algorithm}

\paragraph{Application: Flow around a cylinder.} To demonstrate wtDMD, we consider the classical example of laminar vortex shedding behind a circular cylinder \cite{taira2020modal}.
The cylinder has unit diameter, and the Reynolds number is set to $\mathrm{Re} = 100$,
beyond the critical value where a supercritical Hopf bifurcation produces periodic vortex shedding
\cite{jackson1987finite,zebib1987stability}. At this Reynolds number, the flow settles on a stable limit cycle.
We analyse vorticity data on a rectangular computational domain of size $18 \times 5$,
with no-slip boundary conditions on the cylinder and on the top and bottom walls.
A parabolic inflow is imposed at the left boundary and a non-reflecting outflow at the right.
Each snapshot of vorticity consists of $158,624$ numerical values.

Let $\mathbb{V} \in \mathbb{R}^{158624 \times (N+1)}$ denote the vorticity data.
Because the dynamics are nearly periodic, the data are effectively low-rank.
We project onto a random orthonormal basis $U \in \mathbb{R}^{158624 \times r}$ with target rank $r \geq 1$. The snapshots $X_1,\dots,X_N,X_{N+1}$ are the columns of the matrix $U^\top \mathbb{V} \in \mathbb{R}^{r \times (N+1)}$, from which we form $\mathbb{X},\mathbb{Y}\in\mathbb{R}^{r \times N}$ as in \eqref{SnapshotMatrices}. We use $r \in \{11,21\}$ for illustration, though other choices give similar results.
While $U$ could be taken as the leading left singular vectors of $\mathbb{V}$ (the proper orthogonal modes), we instead fix $U$ to isolate the effects of weighting as $N$ varies.

\begin{figure}[t]  
    \center
    \begin{minipage}[t]{0.015\textwidth}
        \rotatebox[origin=left]{90}{\ \ \ \ \ \ \ \ \ $r = 11$}
    \end{minipage}
    \begin{minipage}[t]{0.32\textwidth}
        \centering
        \includegraphics[width = \textwidth]{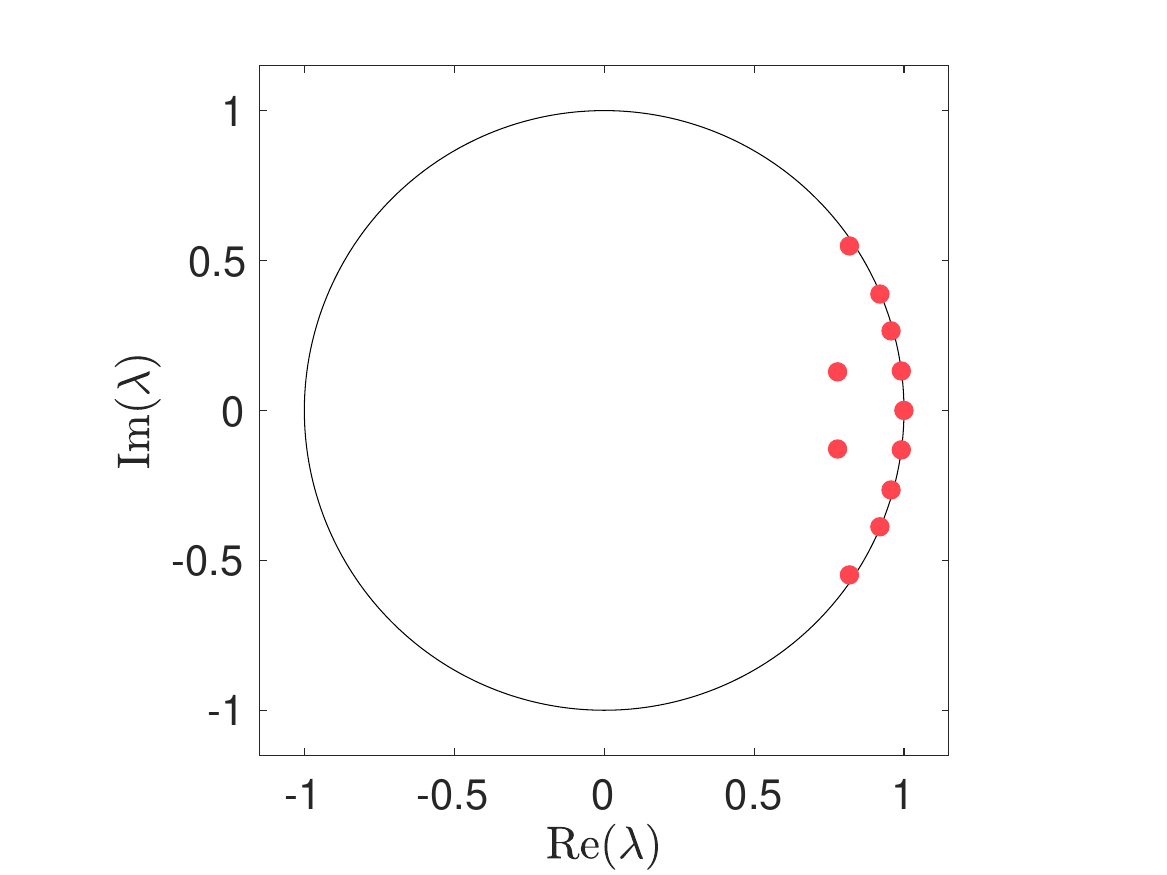}
    \end{minipage}
    \begin{minipage}[t]{0.32\textwidth}
        \centering
        \includegraphics[width = \textwidth]{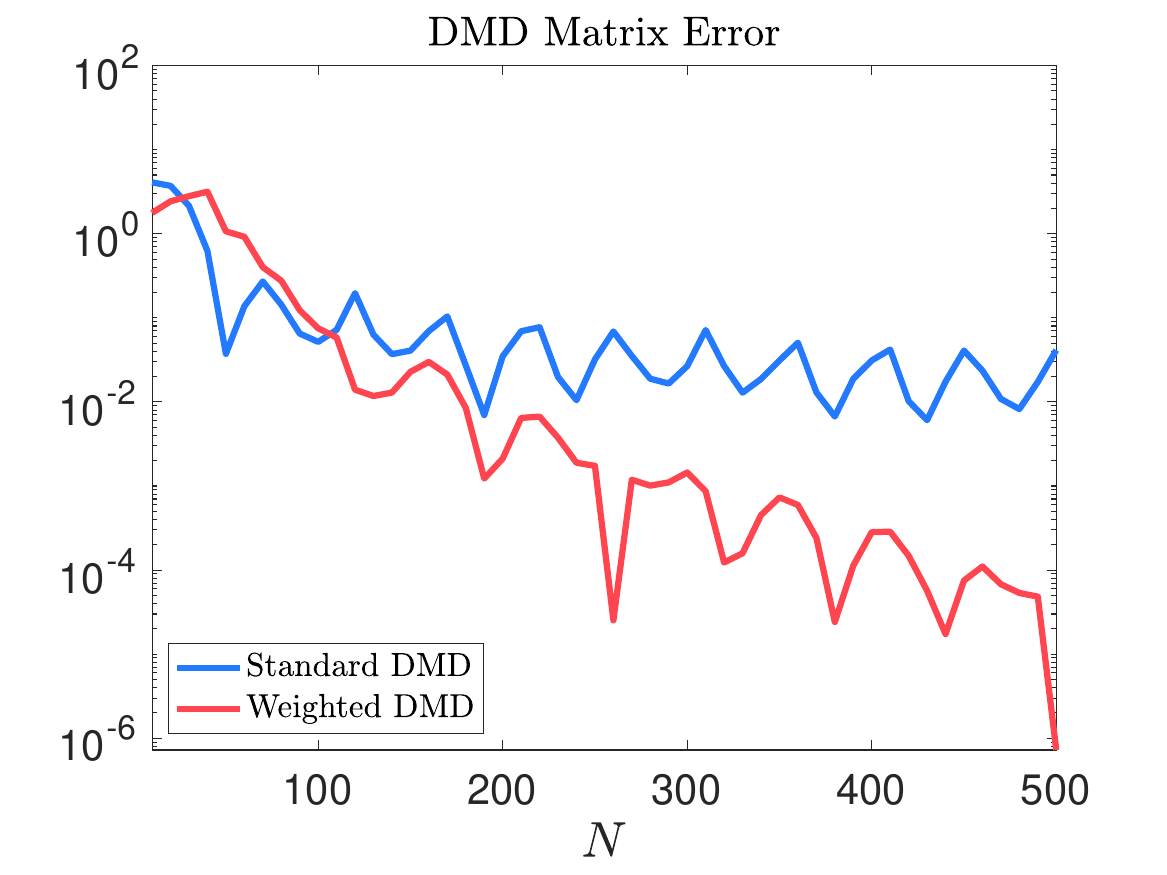}
    \end{minipage}
    \begin{minipage}[t]{0.32\textwidth}
        \centering
        \includegraphics[width = \textwidth]{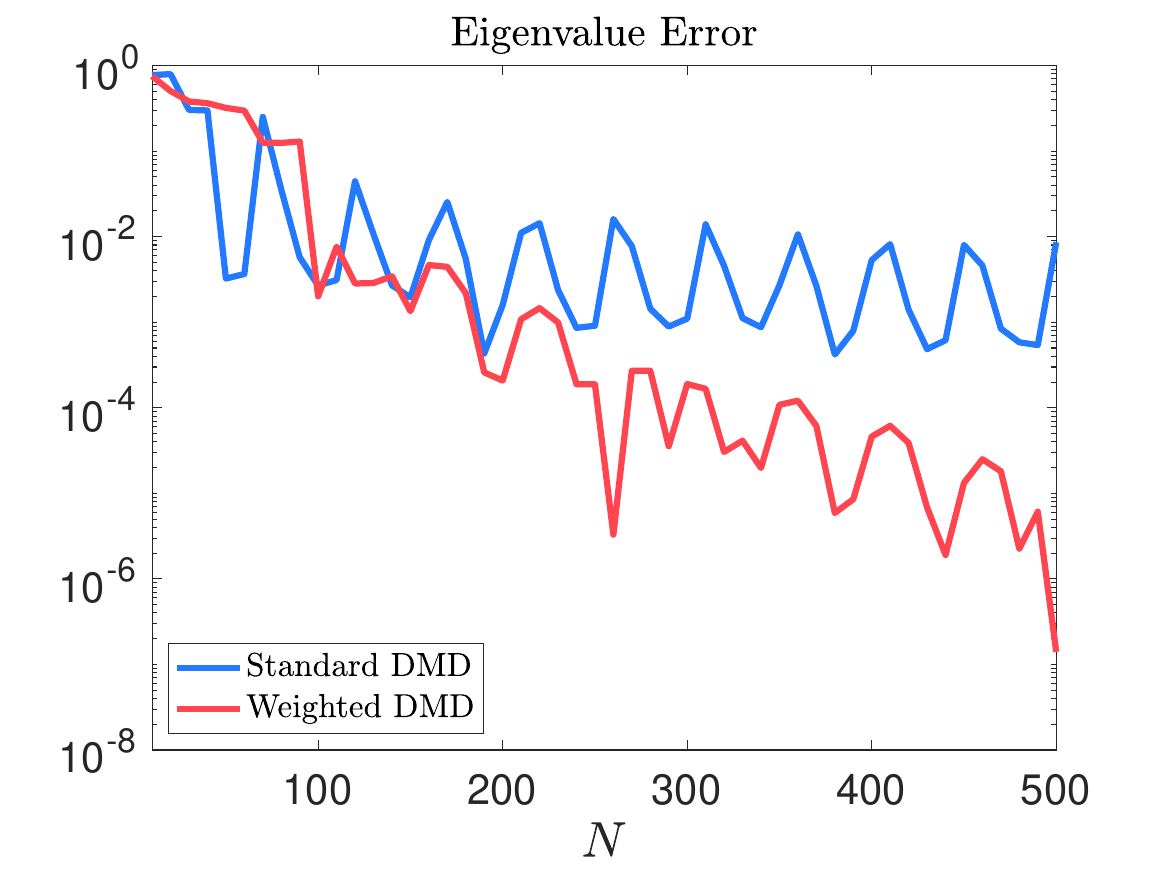}
    \end{minipage}
     \begin{minipage}[t]{0.015\textwidth}
        \rotatebox[origin=left]{90}{\ \ \ \ \ \ \ \ \ \  $r = 21$}
    \end{minipage}
    \begin{minipage}[t]{0.32\textwidth}
        \centering
        \includegraphics[width = \textwidth]{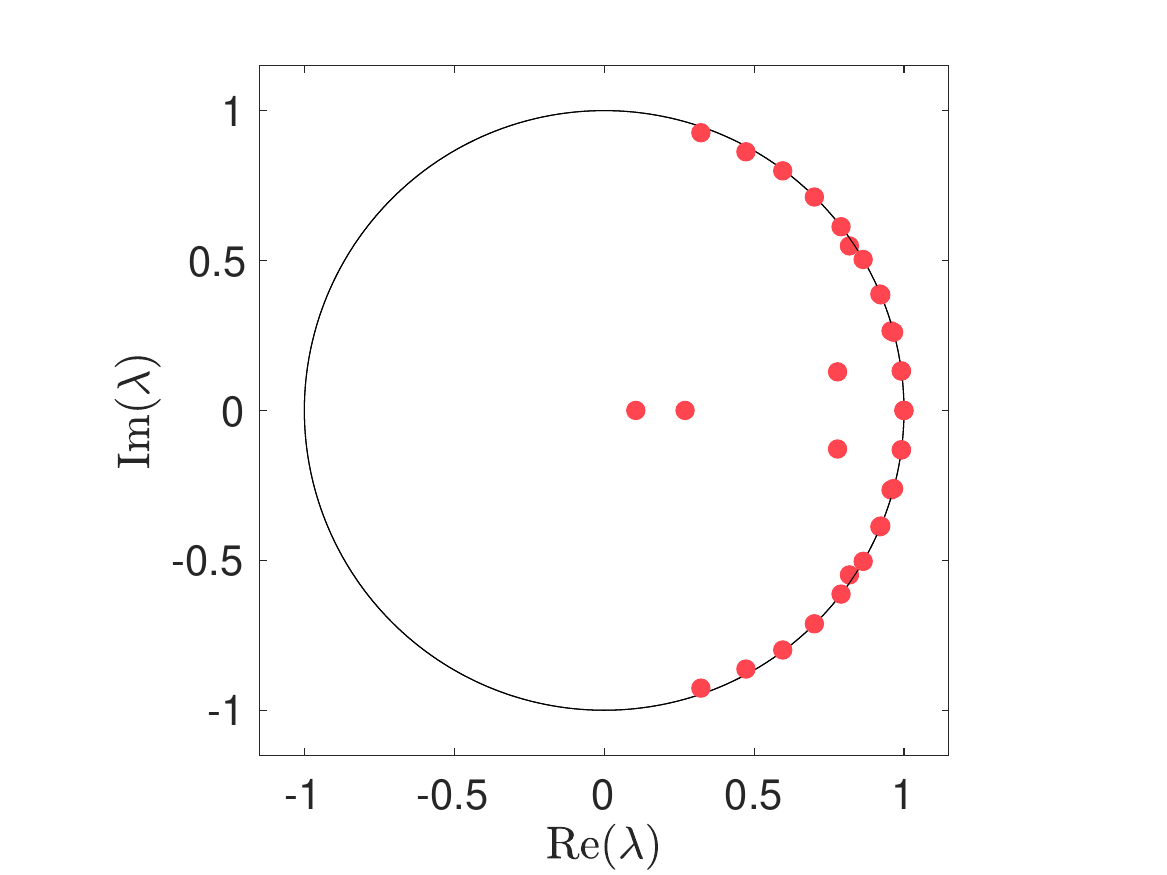}
    \end{minipage}
    \begin{minipage}[t]{0.32\textwidth}
        \centering
        \includegraphics[width = \textwidth]{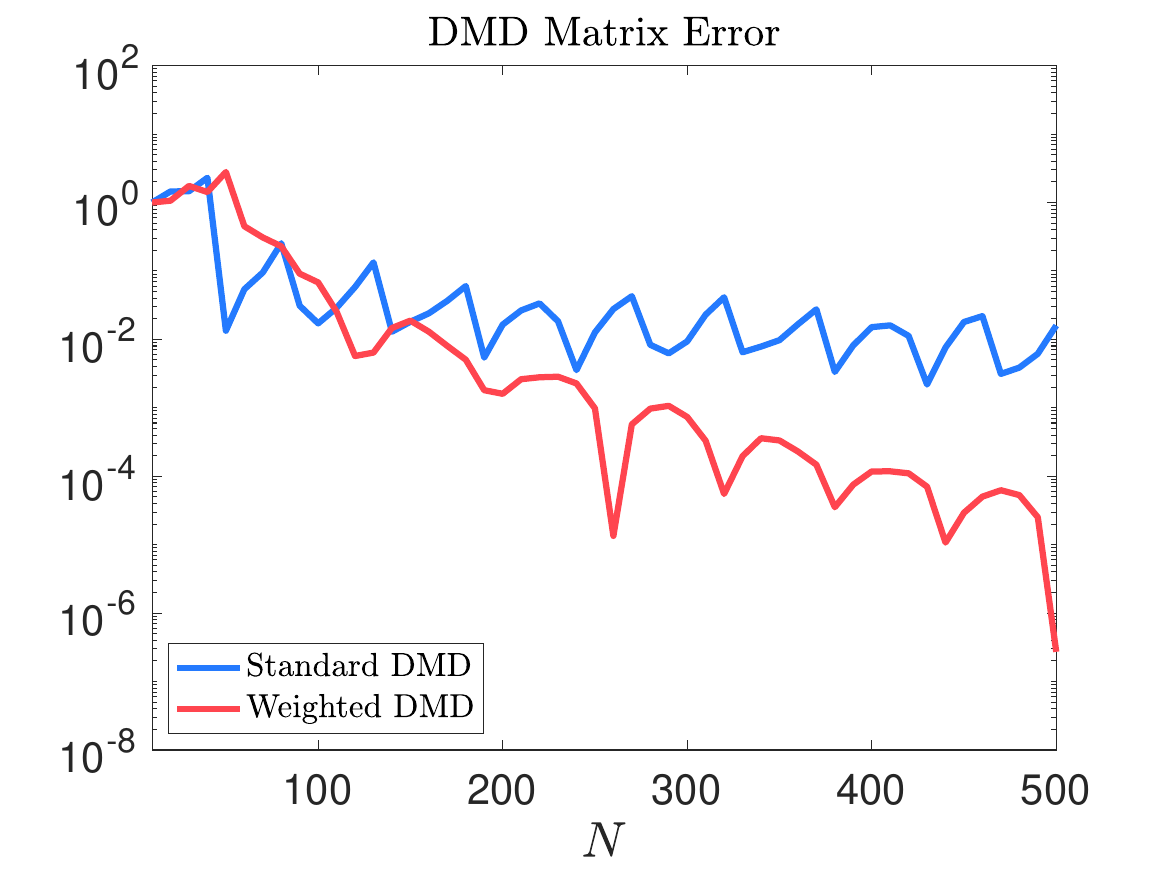}
    \end{minipage}
    \begin{minipage}[t]{0.32\textwidth}
        \centering
        \includegraphics[width = \textwidth]{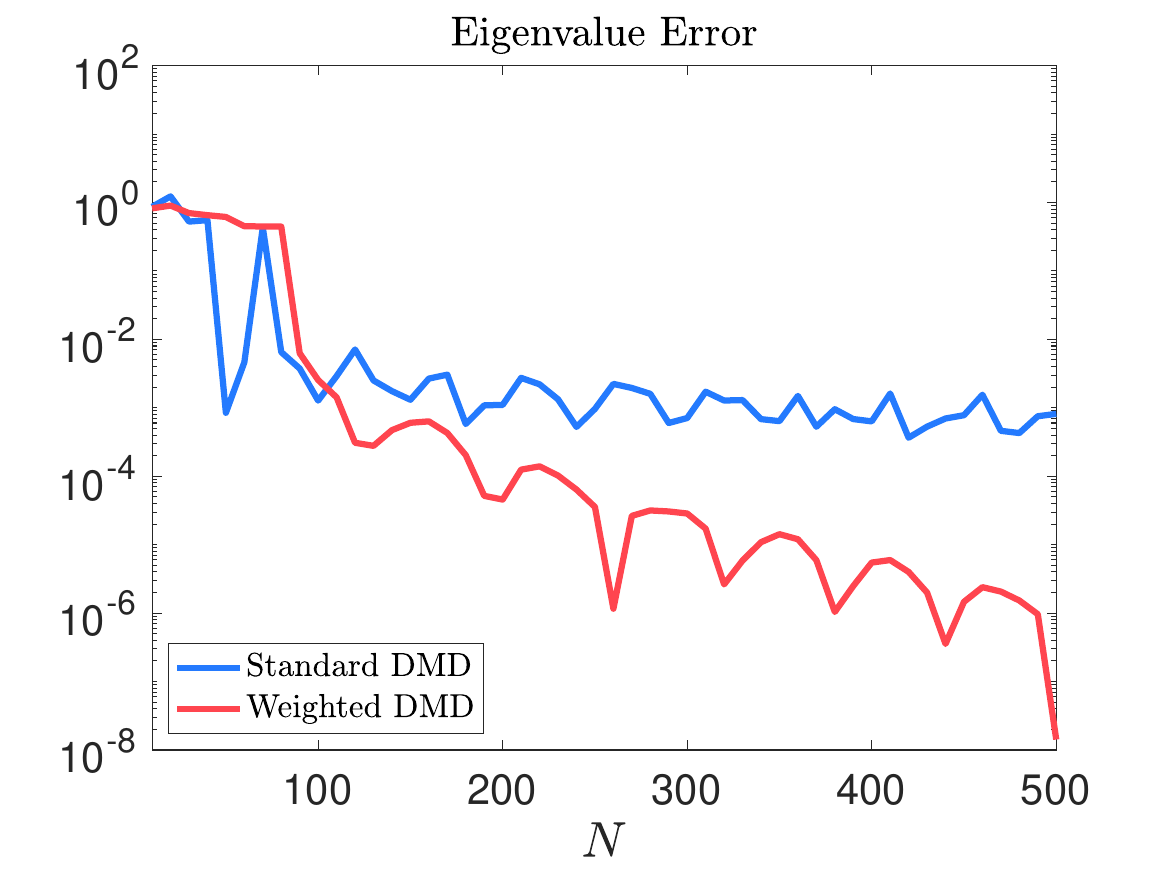}
    \end{minipage}
    \caption{Left: DMD eigenvalues for $N = 1000$ snapshots of vorticity data from flow around a cylinder with $r = 11$ (top) and $r = 21$ (bottom). Middle and Right: Relative errors between the DMD (blue) and wtDMD (red) matrices (middle) and eigenvalues (right) over a range of snapshots $N = 10,20,\dots,500$ and the benchmark results using $N = 1000$ snapshots. The top panels use $r = 11$, while the bottom use $r = 21$.} 
    \label{fig:DMD}
\end{figure}

Taking $N = 1000$, the DMD and wtDMD matrices differ by $\mathcal{O}(10^{-3})$ in Frobenius norm for both $r = 11$ and $r=21$. Their eigenvalues are visually indistinguishable (left column of Figure \ref{fig:DMD}).
A trivial mode at $\lambda=1$ corresponds to the mean flow, while the remaining modes lie approximately on the unit circle, representing oscillatory structures at harmonic multiples of the vortex-shedding frequency.
We use the $N=1000$ results as benchmarks to compare accuracy for smaller $N$.

The middle and right columns of Figure \ref{fig:DMD} show the relative errors of the DMD and wtDMD matrices (middle) and their eigenvalues (right) for $N = 10,20,\ldots,500$, {\color{black} measured against the benchmark $N = 1000$. To compute eigenvalue errors, the $r$ eigenvalues of each DMD matrix are first sorted in decreasing order according to the real-valued quantity $\lambda \mapsto |1-\lambda|$, and the relative error is then taken with respect to the correspondingly sorted eigenvalues for $N=1000$.} For small datasets ($N \lesssim 100$), both methods perform similarly, but once $N \gtrsim 200$, the wtDMD errors drop by several orders of magnitude. Thus, weighting provides faster convergence without altering the asymptotic limit, an easy gain in accuracy from the same data.

\subsection{Example Method II: Approximating the Koopman operator}

Consider the discrete-time stochastic process defined by the random map $x_{n+1} = f(\omega(n),x_n)$, where $\omega\mapsto f(\omega,\cdot)$ is a random variable from a probability space $(\Omega,\mathcal{F},P)$ to the space of continuous maps from $\mathcal{X}\subseteq \mathbb{R}^d$ to itself. The (stochastic) Koopman operator \cite{wanner2022robust,vcrnjaric2020koopman} acts on bounded continuous observables $\phi:\mathcal{X} \to \mathbb{R}$ via
\begin{equation}\label{stochasticKoopman}
    [\mathcal{K}\phi](x) = \int_\Omega\phi(f(\omega,x))\,\mathrm{d}P(\omega) \qquad \forall x \in X.  
\end{equation}
This formulation includes deterministic dynamics $x_{n+1} = f(x_n)$ as a special case. In that setting, the Koopman operator acts by composition \cite{mezic2021koopman,mauroy2020koopman}:
\begin{equation}\label{KoopmanDeterministic}
    [\mathcal{K}\phi](x) = \phi(f(x)) \qquad \forall x \in \mathcal{X}.    
\end{equation}
The Koopman operator provides a linear, deterministic representation of dynamics that may themselves be nonlinear and even stochastic. The price we pay is dimensional: instead of evolving a point in the finite-dimensional state space $\mathcal{X}$, we now evolve observables $\phi:\mathcal{X} \to \mathbb{C}$ in an infinite-dimensional function space. Over the past decade, there has been intense interest in approximating the Koopman operator from data, both for forecasting and for extracting coherent structures; see \cite[Chap.~11]{colbrook_book} and the reviews \cite{brunton2022modern,colbrook2024multiverse}.

The most widely used data-driven method for approximating the Koopman operator on a finite set of observables is extended dynamic mode decomposition (EDMD) \cite{williams2015data}. Suppose we are given sequential data $\{X_n\}_{n = 1}^{N+1}\subset \mathcal{X}$ from a (possibly stochastic) process. We select two dictionaries of bounded, continuous observables,
$$
    \vec\phi := \begin{pmatrix}
       \phi_1 & \cdots & \phi_R 
    \end{pmatrix}, \qquad \vec\psi := \begin{pmatrix}
       \psi_1 & \cdots & \psi_L
    \end{pmatrix},
$$
for some $R,L \geq 1$. EDMD seeks a linear approximation of the projected Koopman operator $\mathcal{P}_L\mathcal{K}\mathcal{P}_R$, where $\mathcal{P}_{L}$ and $\mathcal{P}_{R}$, are the orthogonal projections onto 
$\mathrm{span}\{\vec\psi\}$ and $\mathrm{span}\{\vec\phi\}$, respectively, with orthogonality understood in $L^2(\mu)$, and $\mu$ the ergodic measure of the underlying system. To build this approximation, we form the data matrices
\begin{equation}\label{DictionaryData}
    \Psi = \begin{bmatrix}
        \vec\psi(X_1) \\ \vdots \\ \vec\psi(X_N)
    \end{bmatrix} \in \mathbb{C}^{N \times L}, \qquad \Phi = \begin{bmatrix}
        \vec\phi(X_2) \\ \vdots \\ \vec\phi(X_{N+1})
    \end{bmatrix} \in \mathbb{C}^{N\times R},        
\end{equation}
and set 
\begin{equation}\label{EDMD}
    K^{(N)} := \Psi^\dagger\Phi \in \mathbb{C}^{L\times R}.
\end{equation}
This choice of $K^{(N)}$ minimizes $\|\Phi - \Psi K\|_F$ over all $K \in \mathbb{C}^{L\times R}$. As $N \to \infty$, the matrix $K^{(N)}$ represents the operator $\mathcal{P}_L\mathcal{K}\mathcal{P}_R$ \cite{korda2018convergence,colbrook2023beyond,bramburger2024auxiliary}. When $\vec\phi$ and $\vec\psi$ are just the state coordinates, EDMD reduces (up to transpose) to DMD.
Standard EDMD typically takes $\vec\phi = \vec\psi$ to obtain a square matrix for spectral analysis, but allowing $L \neq R$ can be advantageous; for example, in combining EDMD with semidefinite programming to infer Lyapunov functions and invariant measures directly from data \cite{bramburger2024auxiliary,bramburger2024data}.

\paragraph{Weighted EDMD.}
Just as with DMD, we can rewrite \eqref{EDMD} as
$$
    K^{(N)} = \bigg(\frac{1}{N}\Psi^*\Psi\bigg)^\dagger\bigg(\frac{1}{N}\Psi^*\Phi\bigg),
$$
{\color{black} where the superscript $*$ denotes the conjugate transpose of the matrix,} which exposes the Birkhoff averages that control convergence:
$$
        \frac{1}{N}\Psi^*\Phi = \frac{1}{N}\sum_{n=1}^N \vec\psi(X_n)^*\vec\phi(X_{n+1}), \qquad
        \frac{1}{N}\Psi^*\Psi = \frac{1}{N}\sum_{n=1}^N \vec\psi(X_n)^*\vec\psi(X_n).
$$
This makes it natural to incorporate weighted Birkhoff averages, exactly as in the DMD setting.
We refer to the resulting method as weighted EDMD (wtEDMD).
The goal is the same: accelerate convergence of $\mathcal{K}^{(N)}$ to its asymptotic limit when $N$ is finite.
Algorithm~\ref{alg:wtEDMD} summarises the procedure.

\begin{algorithm}[t]
    \caption{Weighted Extended Dynamic Mode Decomposition (wtEDMD)}\label{alg:wtEDMD}
    {\bf Input:} Sequential snapshots $\{X_n\}_{n = 1}^{N+1}$, bump function $w \in \mathcal{W}$, dictionaries of observables $\vec\phi = (\phi_1, \dots, \phi_R)$, $\vec\psi = (\psi_1, \dots, \psi_L)$.
    \begin{enumerate}
        \item Compute the matrices $\Psi$ and $\Phi$ defined in \eqref{DictionaryData} and $W = \mathrm{diag}(w(0),w(1/N),\dots,w((N-1)/N))$. 
        \item Compute the weighted EDMD matrix $K_w^{(N)} = (W^{1/2}\Psi )^\dagger(W^{1/2}\Phi )$.
    \end{enumerate}
    {\bf Output:} Approximate Koopman operator $\mathrm{span}\{\vec\phi\} \to \mathrm{span}\{\vec\psi\}$ acting as $\mathcal{K}\varphi \approx \vec\psi K_w^{(N)} \vec c $ for all $\varphi =\vec\phi \vec c  \in \mathrm{span}\{\vec\phi\}$ with $\vec c \in \mathbb{C}^R$.
\end{algorithm}

\paragraph{Application to the standard map.} To illustrate wtEDMD, we consider the standard map:
\begin{equation}\label{StandardMap}
        p_{n+1} = p_n + \lambda \sin (\theta_n) \mod 2\pi,\qquad
        \theta_{n+1} = \theta_n +  p_n + \lambda \sin (\theta_n) \mod 2\pi,
\end{equation}
with $\mathcal{X} = \mathbb{T}^2 := \mathbb{S}^1 \times \mathbb{S}^1$ and $\lambda\geq 0$.
We study three regimes: (i) $\lambda = 0.25$, where the trajectory is quasiperiodic; (ii) $\lambda = 5$, where the dynamics are chaotic; (iii) a stochastic case in which $\lambda$ is resampled uniformly from $[0,5]$ at each iteration. We take  $\vec\phi$ and $\vec\psi$ to be Fourier dictionary elements of the form $(p,\theta)\mapsto \mathrm{e}^{\mathrm{i}(k_1 p + k_2 \theta)}$ for $k_1,k_2 \in \{-1,0,1\}$, giving $L = R = 9$. Including more or fewer Fourier modes (changing $L,R$) produces qualitatively similar convergence behaviour.

\begin{figure}[t] 
    \center
    \includegraphics[width = 0.32\textwidth]{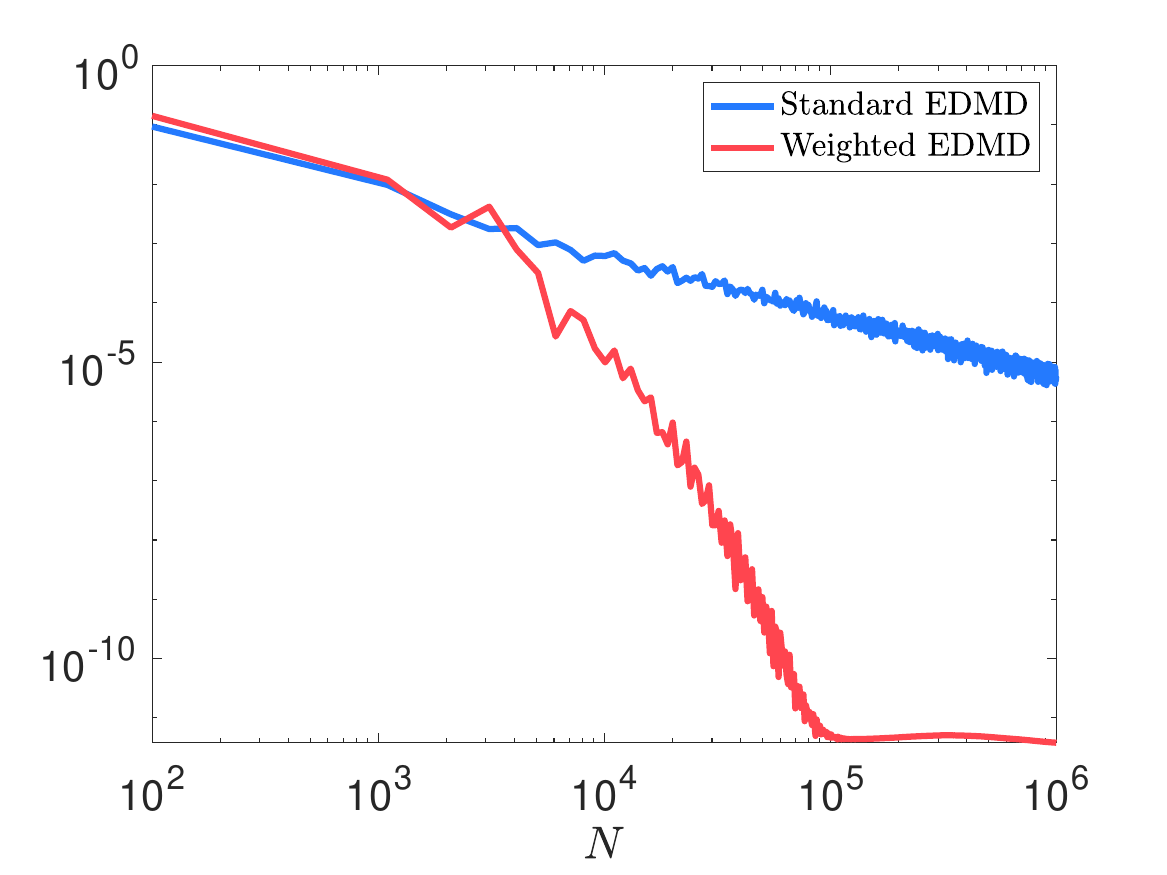}\ 
    \includegraphics[width = 0.32\textwidth]{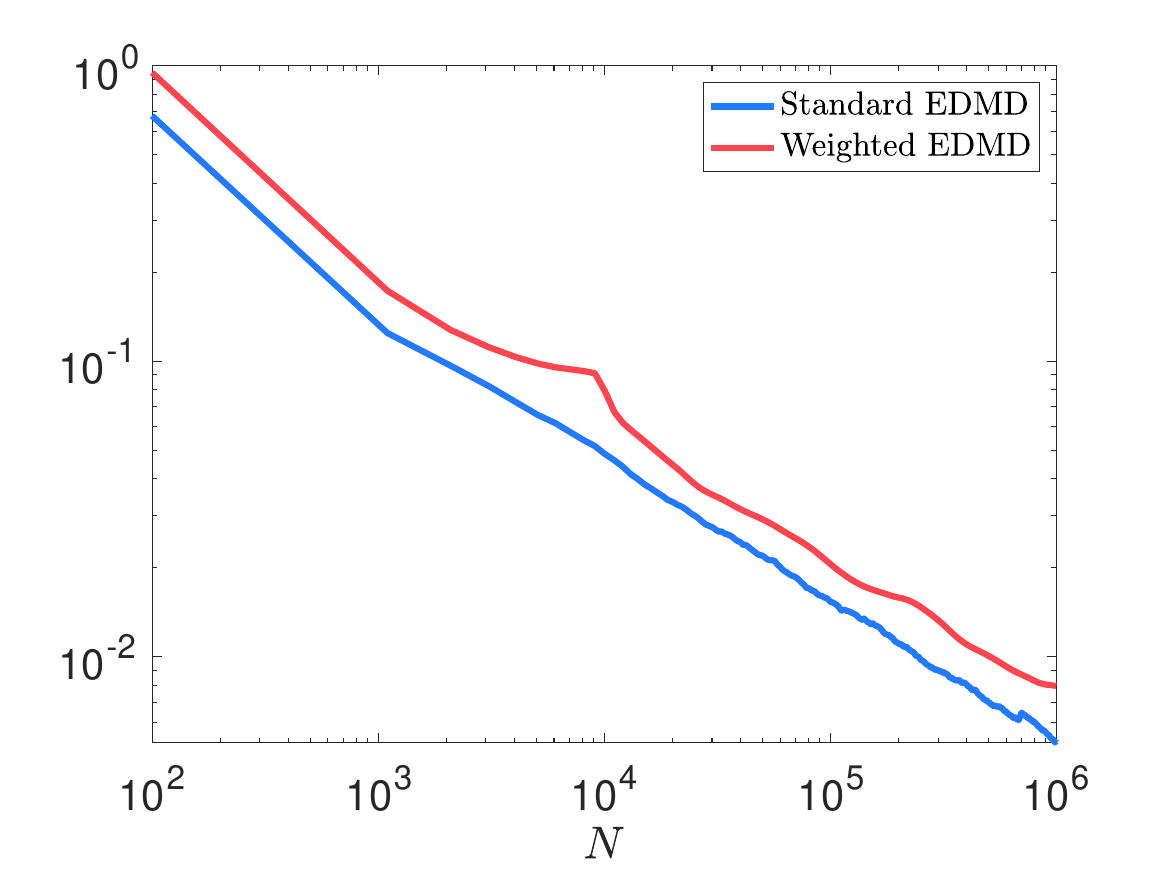} \
    \includegraphics[width = 0.32\textwidth]{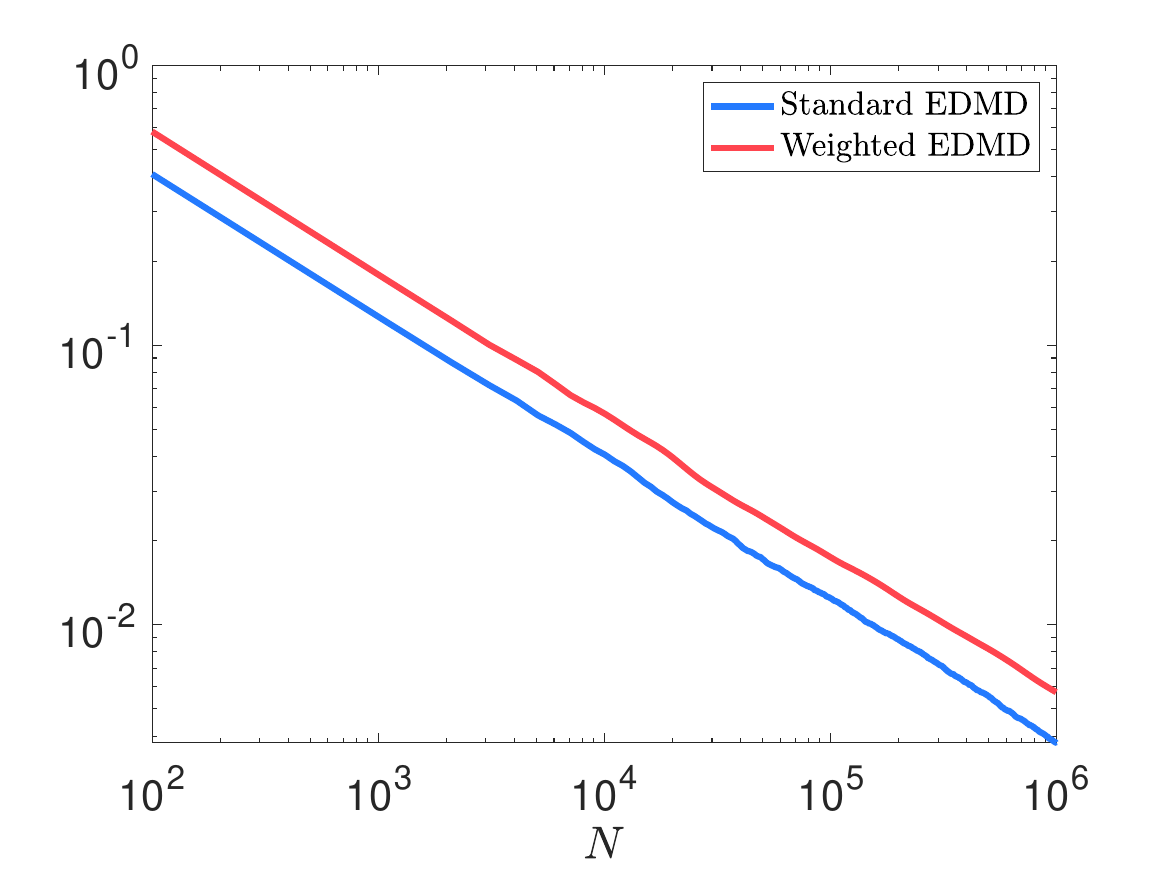} \
    \caption{Averaged relative norm errors between the EDMD (blue) and wtEDMD (red) over a range of snapshots $N$ and the benchmark results using $N = 10^7$ snapshots. All results use data from the standard map \eqref{StandardMap} with: (a) $\lambda = 0.25$; (b) $\lambda = 5$; and (c) $\lambda$ drawn from the uniform distribution on $[0,5]$.}
    \label{fig:EDMD}
\end{figure}

{\color{black} 
We draw a random initial condition $(p_0,\theta_0)$ uniformly from the torus and iterate it under the standard map for $n=10^7$ steps. From these iterations we take the first $N \le 10^6$ snapshots and form the $9\times 9$ EDMD and wtEDMD matrices $K^{(N)}$ and $K_w^{(N)}$. Their relative errors are measured against the $N=10^7$ benchmark by
$$
    \frac{\|K^{(N)} - K^{(10^7)}\|_F}{\|K^{(10^7)}\|_F}, \quad \frac{\|K_w^{(N)} - K_w^{(10^7)}\|_F}{\|K_w^{(10^7)}\|_F},
$$
respectively. We repeat this experiment for 100 independently drawn initial conditions and, for each $N$, average the resulting relative errors. The averaged Frobenius-norm errors are shown in Figure~\ref{fig:EDMD}.} In line with the behaviour seen earlier in Figure~\ref{fig:Logistic}, weighting offers the greatest advantage in the quasiperiodic regime ($\lambda = 0.25$), where convergence is dramatically accelerated. For chaotic and stochastic dynamics, the weighted and unweighted methods converge at essentially the same $\mathcal{O}(1/N)$ rate. Thus, weighting does not degrade performance in difficult regimes, and in favourable regimes it can improve accuracy by orders of magnitude. In applications, replacing EDMD with wtEDMD is therefore a cheap upgrade in accuracy.

\subsection{Example Method III: Model identification}
We now turn to discrete-time model identification.
Given snapshot data $\{X_1,X_2,\dots,X_{N+1}\}\subset\mathbb{C}^d$, we seek a model of the form $x_{n+1} = \Xi\vec\psi(x_n)^\top$, where $\vec\psi$ is a dictionary of candidate nonlinear functions.
This leads to the least-squares problem
\begin{equation}\label{DiscreteModelIdentification}
    \argmin_{\Xi \in \mathbb{C}^{d\times L}} \|\mathbb{Y} - \Xi\Psi^\top\|_F,
\end{equation}
where $\mathbb{Y}$ is as in \eqref{SnapshotMatrices} and $\Psi$ as in \eqref{DictionaryData}. In other words, we are fitting an update rule that advances the state from one snapshot to the next.
As in DMD and EDMD, we can introduce temporal weighting in $\mathbb{Y}$ and $\Psi$ to accelerate convergence and reduce the amount of data required.

\paragraph{SINDy and weighted SINDy.}
The same framework can be adapted to identify continuous-time models.
Suppose we are given snapshots $\{X_1,X_2,\dots,X_{N}\} \subset \mathbb{C}^d$ taken from an underlying ODE $\dot{x}=f(x)$. We approximate the time derivative at each snapshot (e.g. via finite differences) and collect these into $\mathbb{X}'$. We then seek coefficients $\Xi \in \mathbb{C}^{d\times L}$ such that
\begin{equation}\label{ContinuousModelIdentification}
    \argmin_{\Xi \in \mathbb{C}^{d\times L}} \|\mathbb{X}' - \Xi\Psi^\top\|_F  , 
\end{equation}
producing an ODE model of the form $x' = \mathbb{X}'{\Psi^\top}^\dagger\vec\psi(x)^\top$.
The sparse identification of nonlinear dynamics (SINDy) method \cite{brunton2016discovering} adds one more idea: sparsity. Rather than fitting all possible terms in the dictionary, SINDy searches for a model that uses as few terms as possible. Let $\Xi = [\xi_{jk}]_{j = 1,k=1}^{d,L}$. Given a threshold $\eta>0$, small coefficients are pruned and we solve
\begin{equation}\label{SINDy}
    \argmin_{\tilde\Xi= [\tilde\xi_{jk}] \in \mathbb{C}^{d\times L}} \|\mathbb{X}' - \tilde\Xi\Psi^\top\|_F\quad \mathrm{s.t}\quad \tilde\xi_{jk} = 0 \ \mathrm{if}\ |\xi_{jk}| < \eta.
\end{equation}
This procedure is then iterated: update $\xi\leftarrow\tilde\Xi$, re-threshold, and repeat until the coefficients no longer change. Convergence results for SINDy are given in \cite{zhang2019convergence}.

In direct analogy with weighted EDMD, we can apply weights to the data matrices to obtain a weighted SINDy (wtSINDy) method, which aims to produce the same final model but with faster convergence from finite data.
This algorithm is summarised in Algorithm~\ref{alg:wSINDy}.

\begin{algorithm}[t]
    \caption{Weighted Sparse Identification of Nonlinear Dynamics (wtSINDy) for learning continuous time models}\label{alg:wSINDy}
    {\bf Input:} Sequential snapshots $\{(X_n,X_n')\}_{n = 1}^{N}$ of system state and derivative, bump function $w \in \mathcal{W}$, dictionary $\vec\psi = (\psi_1\, \cdots \, \psi_L)$, threshold parameter $\eta > 0$, maximum number of iterates $k_{\mathrm{max}}$.
    \begin{enumerate}
        \item Compute the matrices $\mathbb{X}' = [X_1'\ \dots\ X_N']$, $\Psi$ defined in \eqref{DictionaryData}, and $W = \mathrm{diag}(w(0),w(1/N),\dots,w((N-1)/N))$. 
        \item Initiate $\Xi = (\mathbb{X}'W^{1/2})(\Psi^\top W^{1/2})^\dagger $ and $k=0$.
        \item  Perform sequential least-squares:
        \begin{enumerate}
            \item Compute set of indices $\mathcal{I}=\{(j,k):|\Xi_{jk}|\geq\eta\}$.
            \item Set $\Xi_{\mathrm{new}}$ to be the $d\times L$ zero matrix.
            \item For $j=1,\ldots,d$, set $\mathcal{I}_j=\{k:(j,k)\in\mathcal{I}\}$ and update
            $$
            \Xi_{\mathrm{new}}(j,\mathcal{I}_j)\leftarrow
            [\mathbb{X}'W^{1/2}](j,:)([\Psi^\top W^{1/2}](\mathcal{I}_j,:))^\dagger.
            $$
            \item If $\Xi_{\mathrm{new}}=\Xi$ or $k\geq k_{\mathrm{max}}$ terminate, otherwise update $\Xi\leftarrow \Xi_{\mathrm{new}}$ and $k\leftarrow k+1$, and return to step (a).
        \end{enumerate}
    \end{enumerate}
    {\bf Output:} Coefficient matrix $\Xi$ and sparse ODE model $x' = \Xi\vec\psi(x)^\top$.
\end{algorithm}

\begin{figure}[t] 
    \center
    \includegraphics[width = 0.49\textwidth]{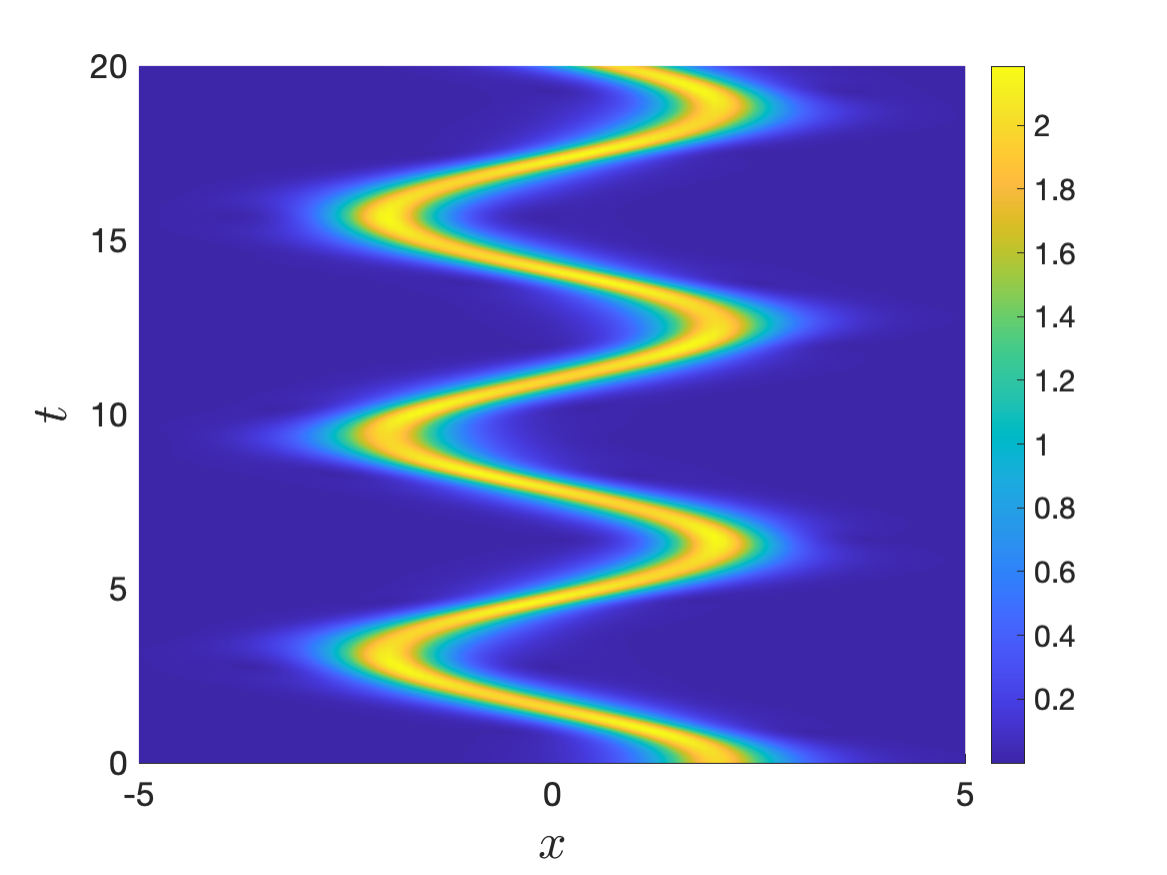}\hfill
    \includegraphics[width = 0.49\textwidth]{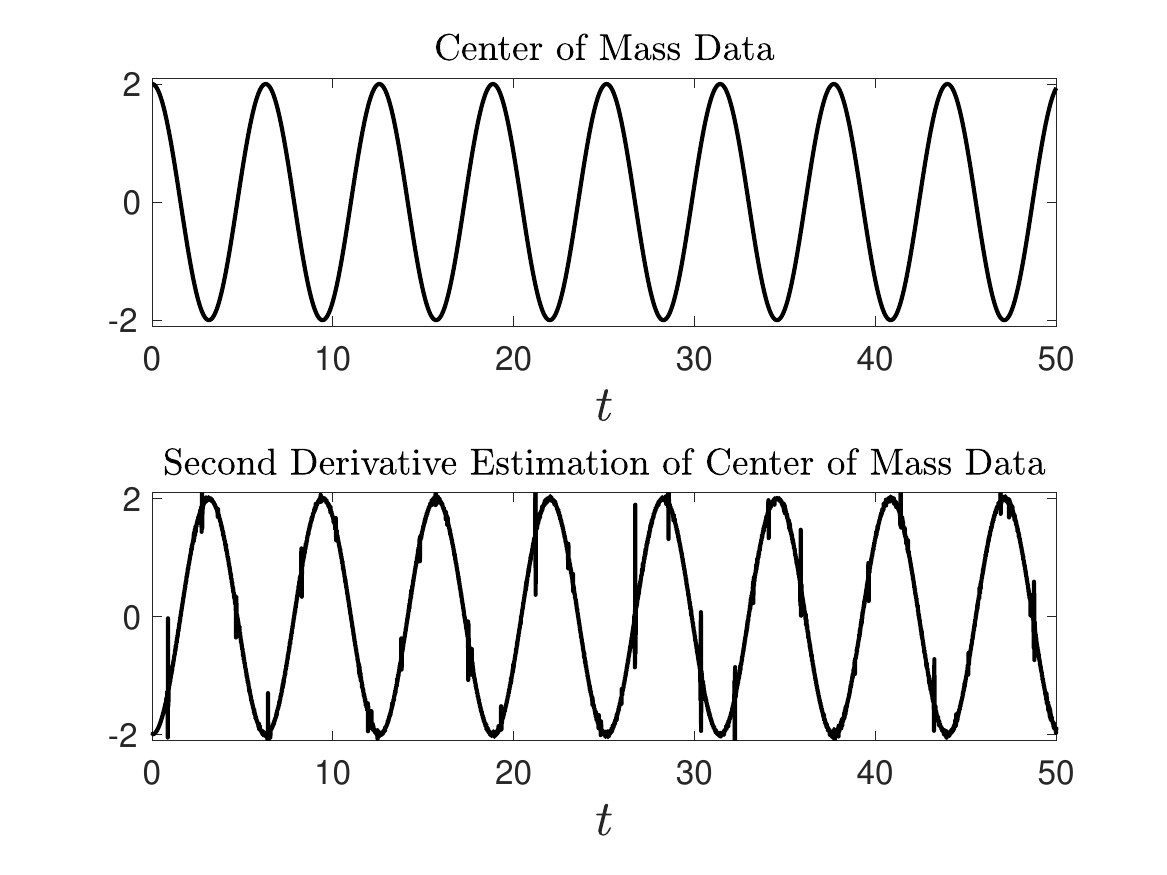}
    \caption{Left: A top-down surface plot of a solution to the NLSE \eqref{NLSE} initialised with $u(x,0) = 2\,\mathrm{sech}(2 (x - x_0))$. Right: Model identification is performed on the tracked centre of mass of the soliton on the left (top) and numerical estimations of its second derivative (bottom).}
    \label{fig:NLSE}
\end{figure}

\paragraph{Application: Harmonic oscillations in the NLSE.} We now test wtSINDy on a reduced model of soliton motion in the nonlinear Schrödinger equation (NLSE) with a parabolic trap,
\begin{equation}\label{NLSE}
    \mathrm{i}\frac{\partial u}{\partial t} = -\frac{1}{2}\frac{\partial^2 u}{\partial x^2} - |u|^2u + \frac{x^2}{2}u, 
\end{equation}
where $u = u(x,t)$ is complex-valued, $t \geq 0$, $x \in [-5,5]$, and periodic boundary conditions are imposed.
As discussed in \cite{kevrekidis2005statics}, if we initialise with $u(x,0) = 2\mathrm{sech}(2 (x - x_0))$ for some $x_0 \neq 0$, the resulting soliton oscillates harmonically with frequency 1. This behaviour follows from the Ehrenfest theorem \cite{hall2013quantum}. The left panel of Figure~\ref{fig:NLSE} shows $|u(x,t)|$ for $x_0 = 2$, exhibiting this oscillation.

We want to identify a low-dimensional model that captures this motion.
To do so, we track the soliton's centre of mass:
\begin{equation}
    X_n = \frac{\int_{-5}^{5} x|u(x,nk)|\,\mathrm{d}x}{\int_{-5}^{5} |u(x,nk)|\,\mathrm{d}x}, \quad n = 0,\dots, N+1,
\end{equation}
where the NLSE \eqref{NLSE} is integrated spectrally, the spatial integrals are discretised on 128 uniformly spaced points in $[-5,5]$, and data are sampled at time step $k = 0.01$. We aim to learn a second-order model of the form $x'' = f(x)$,  whose ground truth is $f(x) = -x$, i.e. the harmonic oscillator

To approximate the second derivative, we compute
\begin{equation}\label{DerivativeEstimate}
    X_n'' := \frac{X_{n+1} + X_{n-1} - 2X_n}{k^2}, \quad n = 1,\dots, N,
\end{equation}
and take as our dictionary $\vec\psi = \{1,x,x^2,x^3,x^4,x^5\}$. (The results are insensitive to the precise degree.) The right panel of Figure~\ref{fig:NLSE} shows the resulting position data and noisy second-derivative estimates. Note in particular that estimating the second derivative by finite differences introduces visible noise. 

\begin{figure}[t] 
    \center
    \includegraphics[width = 0.43\textwidth]{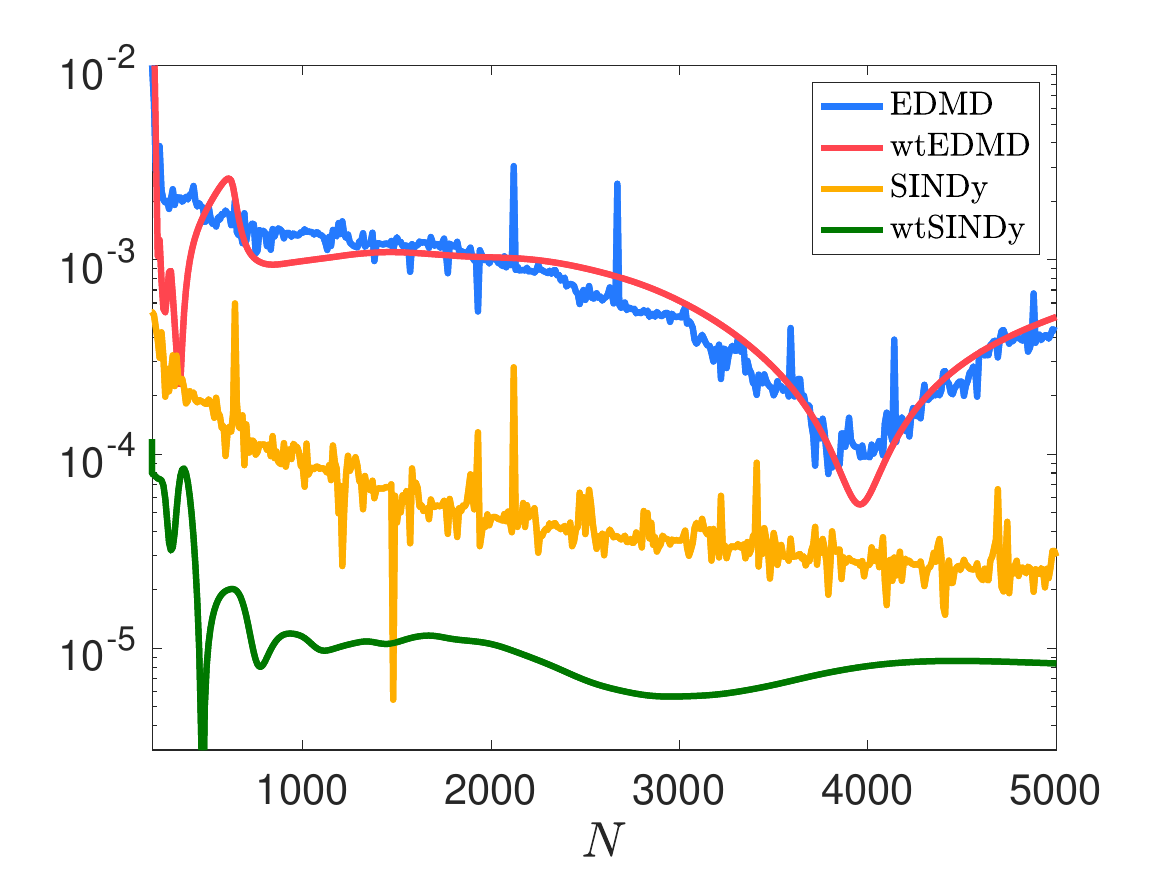}\hfill
    \includegraphics[width = 0.43\textwidth]{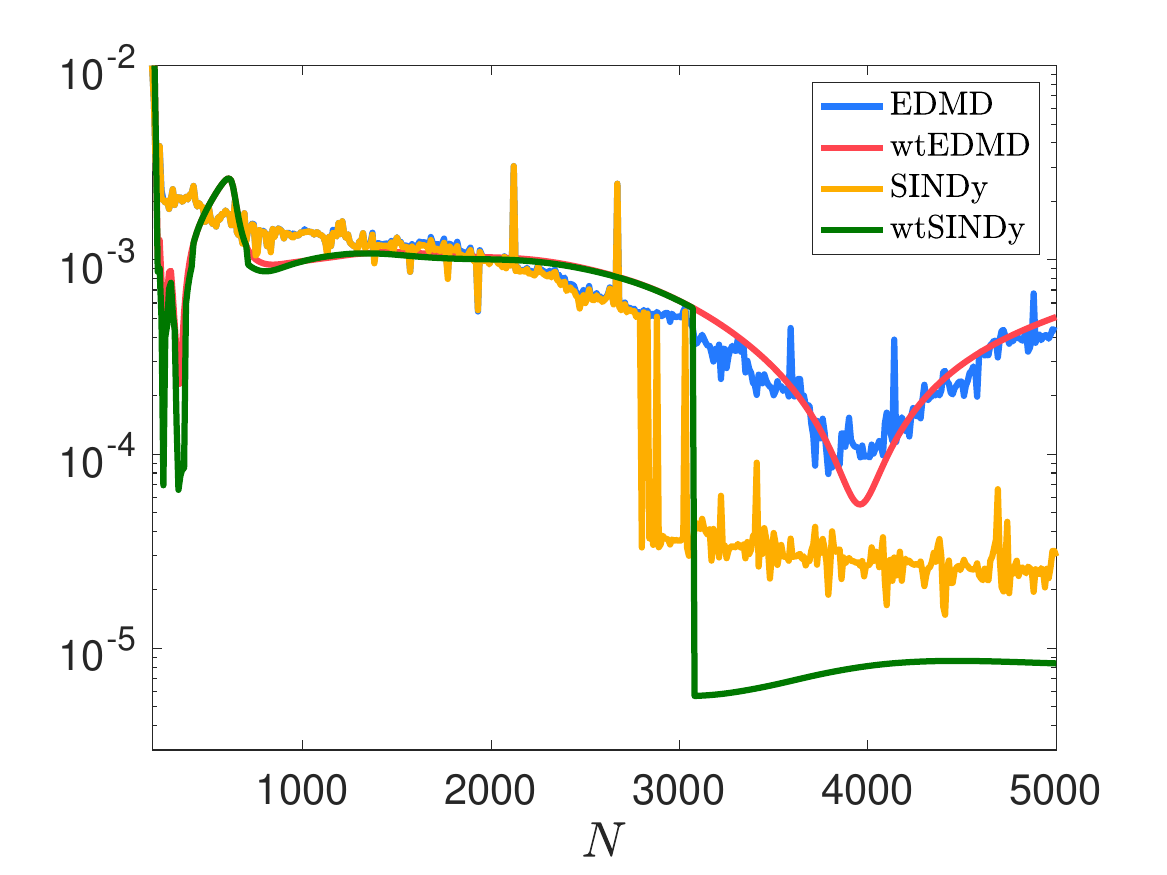}
    \caption{Relative error over $N\geq 200$ snapshots for different model identification algorithms to discover the harmonic oscillator $x'' = -x$ using the centre of mass data shown in Figure~\ref{fig:NLSE}. SINDy (yellow) and wtSINDy (green) methods use a sparsity parameter of $\eta = 10^{-2}$ (left) and $\eta = 10^{-4}$ (right). 
    }
    \label{fig:wSINDy}
\end{figure}

Figure~\ref{fig:wSINDy} reports the error in the recovered model. We compare: least-squares \eqref{ContinuousModelIdentification} (labelled EDMD), its weighted variant (wtEDMD), SINDy \eqref{SINDy}, and weighted SINDy (wtSINDy). {\color{black} For each method, we measure the Frobenius-norm error $\|\Xi^{(N)} - \Xi_\mathrm{exact}\|_F$ between the coefficient vector $\Xi^{(N)}$ identified from $N$ snapshots and the exact coefficients $\Xi_\mathrm{exact} = [ 0 \, {-}1 \, 0 \, 0 \, 0 \, 0 ]$.} The left panel uses a sparsity threshold $\eta = 10^{-2}$; the right panel uses $\eta = 10^{-4}$.

EDMD and wtEDMD perform similarly, with no clear gain or loss from weighting.
Although the motion is periodic, the limited benefit of weighting here is likely due to numerical error in computing the centre of mass, and especially the second derivatives.
Both SINDy methods outperform EDMD and wtEDMD because the sparsification step suppresses this noise.
SINDy and wtSINDy agree until the coefficients become small enough to be thresholded.
Beyond that point, wtSINDy achieves errors an order of magnitude smaller than SINDy.

\subsection{Example Method IV: Spectral measures}

One advantage of viewing a dynamical system through the Koopman operator is that linearity opens the door to spectral analysis of global behaviour.
For clarity, consider the deterministic system $x_{n+1} = f(x_n)$, so that the Koopman operator acts by composition as in \eqref{KoopmanDeterministic}. If $\phi_\lambda$ is an eigenfunction of $\mathcal{K}$ with eigenvalue $\lambda \in\mathbb{C}$, then
\begin{equation}\label{KoopmanCoherent}
    \mathcal{K}\phi_\lambda(x_n) = [\mathcal{K}^n\phi_\lambda](x_0) = \lambda^n\phi_\lambda(x_0), \quad n = 1,2,3,\dots.
\end{equation}
Thus $\phi_\lambda$ evolves along trajectories by pure exponential growth, decay, or oscillation.
Koopman eigenfunctions reveal coherent structures in the dynamics: their level sets can correspond to invariant manifolds \cite{mezic2015applications}, isostables \cite{mauroy2013isostables}, and other dynamically relevant sets.

To speak precisely about the spectrum, we must specify the space of observables.
Assume that the map $f$ preserves a measure $\mu$ (e.g., an ergodic measure), and consider the Hilbert space $L^2(\mu)$. Then $\mathcal{K}$ is an isometry on $L^2(\mu)$; that is $\|\mathcal{K}\phi\|_{L^2(\mu)} = \|\phi\|_{L^2(\mu)}$ for all observables $\phi \in L^2(\mu)$. If we additionally assume that the dynamics are invertible modulo $\mu$-null sets, then $\mathcal{K}$ is unitary on $L^2(\mu)$. In that case, its spectrum $\sigma(\mathcal{K})$ lies on the unit circle $\mathbb{S}^1 =  \{\lambda \in \mathbb{C}:\ |\lambda| = 1\}$, and the associated eigenfunctions have neutral time evolution: no growth, no decay.

In general, there is no reason for $\mathcal{K}$ to have a complete basis of eigenfunctions.
Nevertheless, the spectral theorem \cite[Thm.~X.4.11]{conway2019course} gives a diagonalization of $\mathcal{K}$
using a projection-valued measure $\mathcal{E}$ supported on $\sigma(\mathcal{K}) \subset \mathbb{S}^1$. This assigns to each Borel subset of $\mathbb{S}^1$ an orthogonal projector and yields
\begin{equation}
    \phi = \bigg(\int_{\mathbb{S}^1} \mathrm{d}\mathcal{E}(\lambda)\bigg)\phi \quad \mathrm{and} \quad \mathcal{K}\phi = \bigg(\int_{\mathbb{S}^1} \lambda \mathrm{d}\mathcal{E}(\lambda)\bigg)\phi,
\end{equation}
for all $\phi \in L^2(\mu)$. Iterating gives the analog of \eqref{KoopmanCoherent}: 
\begin{equation}
    \phi(x_n) = [\mathcal{K}^n\phi](x_0) = \bigg(\int_{\mathbb{S}^1}\lambda^n \mathrm{d}\mathcal{E}(\lambda)\bigg)\phi(x_0) \quad \forall \phi \in L^2(\mu).
\end{equation}
The measure $\mathcal{E}$ decomposes the state space into invariant subspaces indexed by subsets of the spectrum, even when actual eigenfunctions are scarce. For a review of the data-driven approximation of Koopman spectral measures, see \cite[Sec.~6]{colbrook2025introductory}.

Now fix a unit-norm observable $\phi \in L^2(\mu)$, $\|\phi\|_{L^2(\mu)} = 1$.
The associated spectral measure of $\mathcal{K}$ with respect to $\phi$ is the probability measure $\nu_\phi(A) = \langle \mathcal{E}(A)\phi,\phi\rangle$, supported on $\mathbb{S}^1$. 
Spectral measures carry information not only about Koopman eigenfunctions but also about generalsed eigenfunctions \cite{colbrook2024rigorous,colbrook2025rigged}, which in turn illuminate features such as global stability of equilibria \cite{mauroy2016global} and ergodic partitions \cite{budivsic2012applied}.
These objects play a central role in nonlinear dynamics \cite{mezic2021koopman}.

Since $\sigma(\mathcal{K})\subset \mathbb{S}^1$, write $\lambda = \mathrm{e}^{\mathrm{i}\theta}$ and transfer $\nu_\phi$ to a measure $\xi_\phi$ on the periodic interval $[-\pi,\pi]$ so that $\xi_\phi(B) = \nu_\phi(\mathrm{e}^{\mathrm{i}B})$ for measurable subsets $B \subseteq [-\pi,\pi]$. A direct computation using the projection-valued measure gives its Fourier coefficients \cite{colbrook2025rigged}:
\begin{equation}
    \hat\xi_\phi(n) = \overline{\hat\xi_\phi(-n)} = \frac{1}{2\pi}\langle\phi,\mathcal{K}^{n}\phi \rangle, \quad n \in\mathbb{Z}.
\end{equation}
Thus, if we can approximate the autocorrelations $\langle\phi,\mathcal{K}^{n}\phi \rangle$, we can approximate the spectral measure. Given trajectory data $\{X_n\}_{n = 1}^{N}$, we estimate these autocorrelations by the Birkhoff averages
\begin{equation}\label{autocorrelation}
    \langle\phi,\mathcal{K}^{n}\phi \rangle \approx \frac{1}{N-n}\sum_{j = 1}^{N-n}\phi(X_j)\overline{\phi(X_{j + n})} =: a_n.
\end{equation}

\paragraph{Weighted spectral measure approximation.}
These averages can be weighted, just as before, to accelerate convergence and thereby obtain accurate spectral information from shorter trajectories.
To turn the autocorrelations into an approximation of the spectral density, we introduce a filter $\varphi:[-1,1]\to \mathbb{R}$ and use the first $M$ autocorrelations to define
\begin{equation}\label{spectraldensity}
    \xi_{\phi,M}(\theta) = \sum_{n = -M}^M \varphi\bigg(\frac{n}{M}\bigg)a_n\mathrm{e}^{\mathrm{i}n\theta}.
\end{equation}
The filter $\varphi$ should be close to 1 when $x$ is close to $0$ and tapers to $0$ near $x = \pm 1$; see \cite{colbrook2024rigorous} for details.
The full weighted procedure is summarised in Algorithm~\ref{alg:SpecMeas}.

\begin{algorithm}[t]
    \caption{Weighted approximation of spectral measures from data}\label{alg:SpecMeas}
    {\bf Input:} Sequential snapshots $\{X_n\}_{n = 1}^{N+1}$, bump function $w \in \mathcal{W}$, observable $\phi$, Fourier filter $\varphi$, and $M \geq 1$.
    \begin{enumerate}
        \item Approximate the first $M$ autocorrelations \eqref{autocorrelation} using the weighted average $a_n = \frac{1}{2\pi\alpha_{N-n}}\sum_{j = 1}^{N - n} w((j-1)/(N-n))\phi(X_j)\overline{\phi(X_{j + n})}$ for $0 \leq n \leq M$ with $\alpha_{N-n} := \sum_{j = 1}^{N-n} w((j-1)/(N-n))$. 
        \item Set $a_{-n} = \overline{a_n}$ for $1 \leq n \leq M$.
    \end{enumerate}
    {\bf Output:} The density $\xi_{g,M}(\theta) = \sum_{n = -M}^M \varphi(\frac{n}{M})a_n\mathrm{e}^{\mathrm{i}n\theta}$ that can be evaluated at any $\theta \in [-\pi,\pi]$.
\end{algorithm}

\paragraph{Application: Lid-driven cavity flow.}
We now apply this method to two-dimensional lid-driven cavity flow at various Reynolds numbers.
This flow describes an incompressible viscous fluid in a rectangular cavity with a moving lid.
As the Reynolds number increases, the flow undergoes a cascade of bifurcations, exhibiting periodic, quasiperiodic, skew-periodic,\footnote{Here, ``skew-periodic'' denotes a flow with a mixture of (quasi)periodic and chaotic components.} and eventually chaotic/turbulent dynamics \cite{shen1991hopf,arbabi2017study}. We use the dataset from \cite{arbabi2017study}, which provides $20,000$ snapshots of the vorticity field, sampled at time step $0.1$, for Reynolds numbers 
$13,000$ (periodic), 
$16,000$ (quasiperiodic), and 
$19,000$ (skew-periodic).
As our observable $\phi$, we take the total kinetic energy of the flow.

\begin{figure}[t] 
    \center
    \begin{minipage}[t]{0.32\textwidth}
        \centering
        \includegraphics[width = \textwidth]{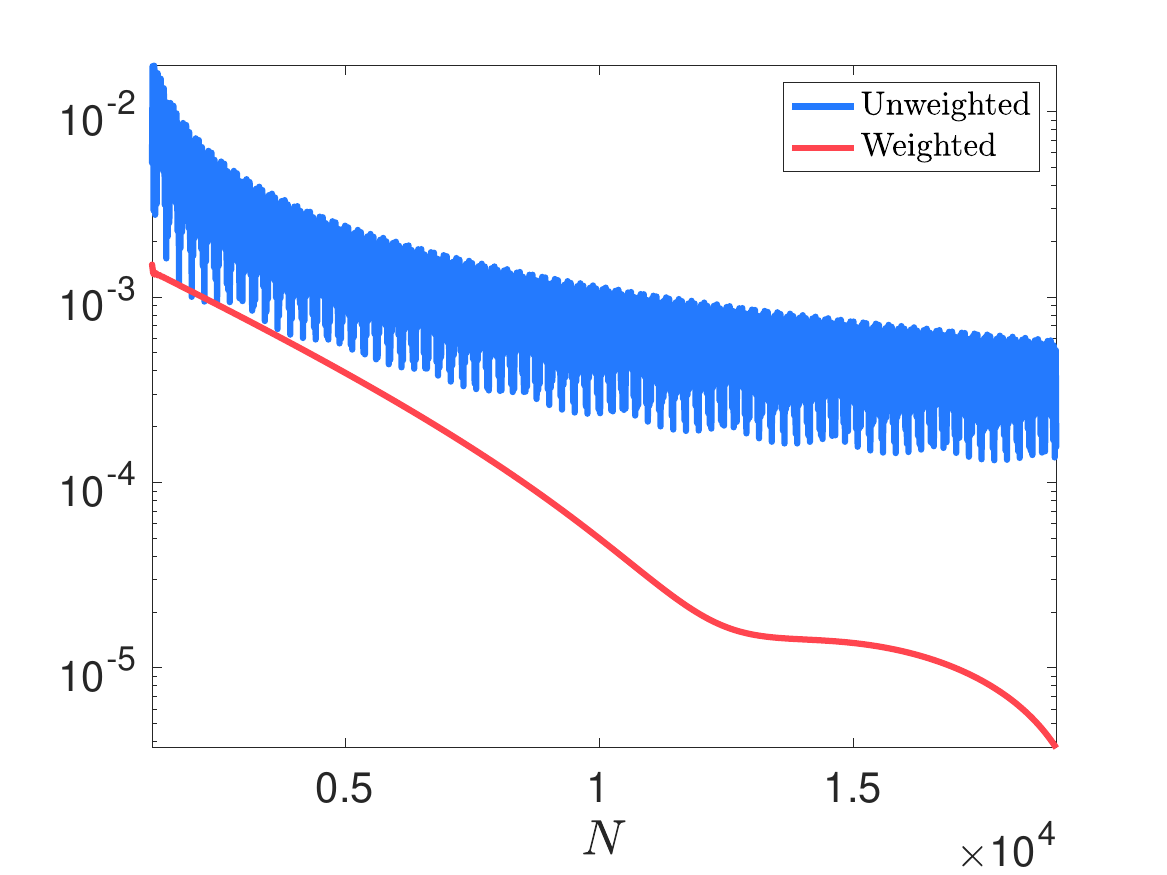}
    \end{minipage}
    \begin{minipage}[t]{0.32\textwidth}
        \centering
        \includegraphics[width = \textwidth]{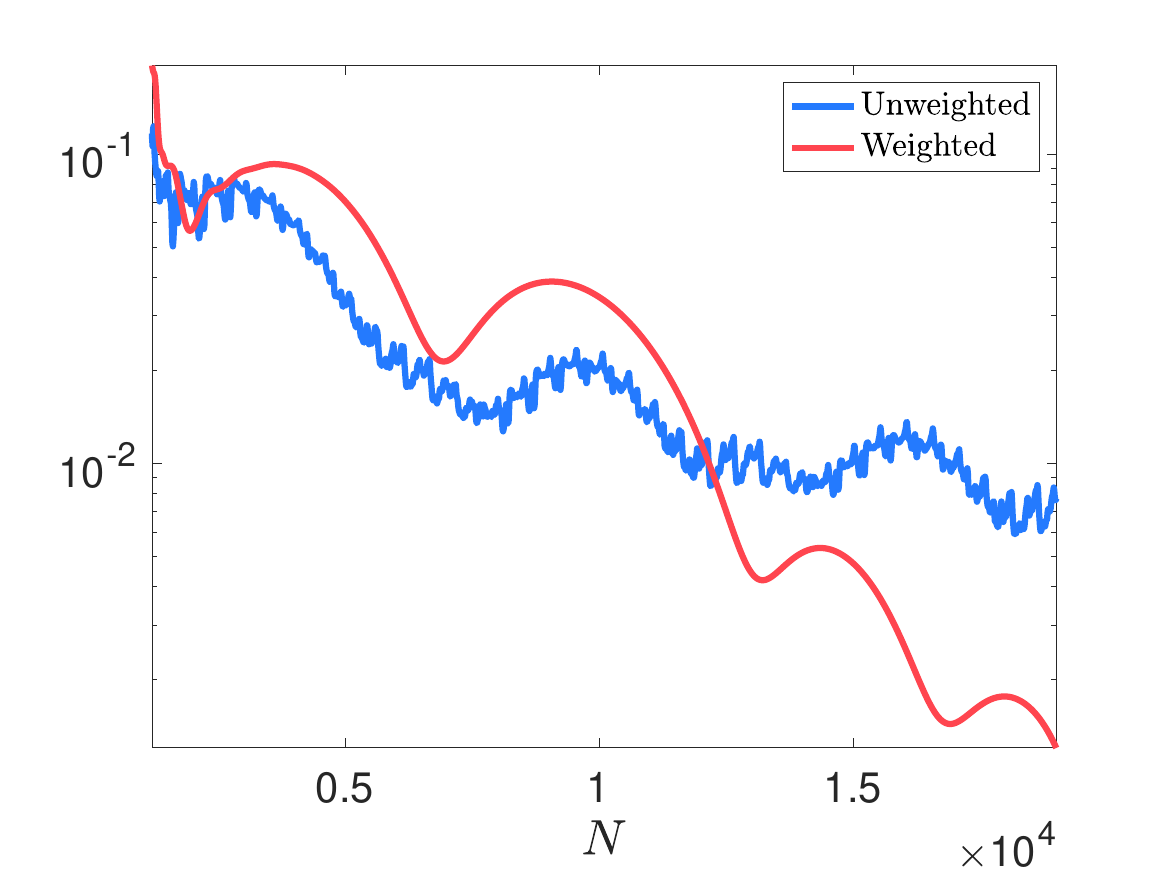}
    \end{minipage}
    \begin{minipage}[t]{0.32\textwidth}
        \centering
        \includegraphics[width = \textwidth]{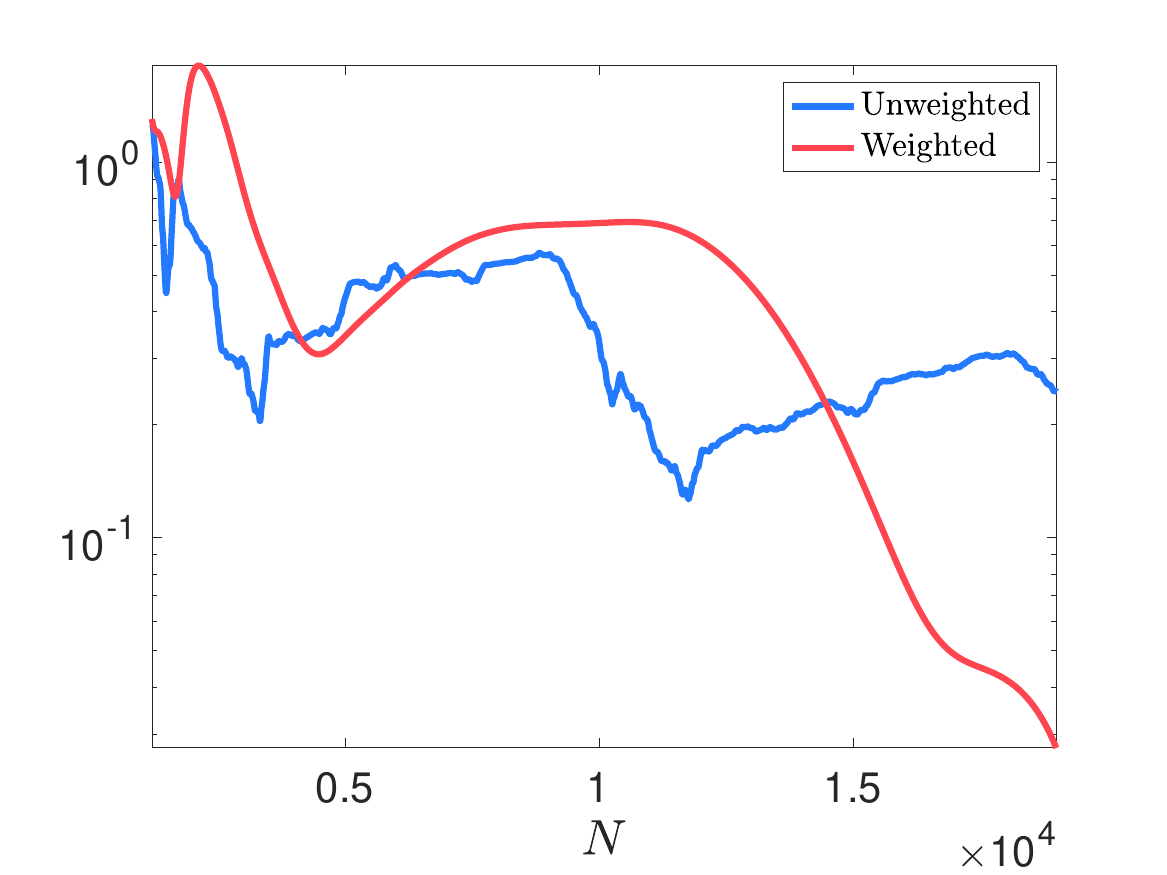}
    \end{minipage}
    \caption{Relative norm errors between the unweighted (blue) and weighted (red) the first $M = 1000$ autocorrelation approximations \eqref{autocorrelation} using the kinetic energy of lid-driven cavity flow over a range of snapshots $N$ and the benchmark results using $N = 20,000$ snapshots. Reynolds numbers for the fluid flow are $13,000$ (left), $16,000$ (middle), and $19,000$ (right).
    }
    \label{fig:auto_error}
\end{figure}

\begin{figure}[t] 
    \center
    \begin{minipage}[t]{0.015\textwidth}
        \rotatebox[origin=left]{90}{\ Unweighted}
    \end{minipage}
    \begin{minipage}[t]{0.32\textwidth}
        \centering
        \includegraphics[width = \textwidth]{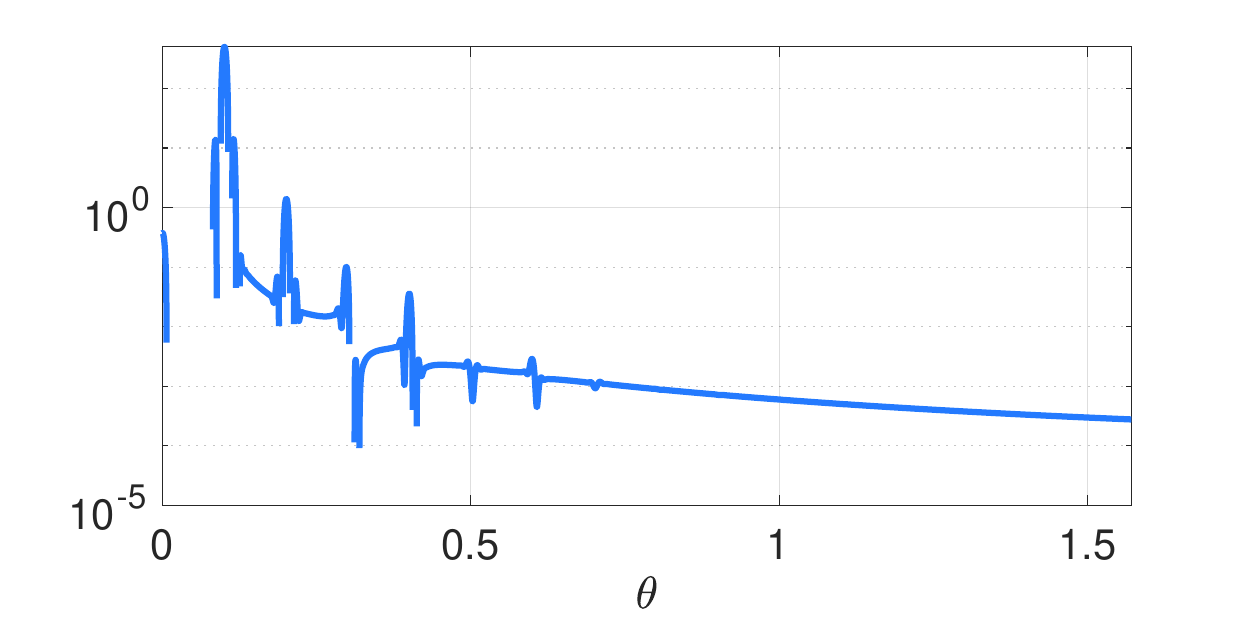}
    \end{minipage}
    \begin{minipage}[t]{0.32\textwidth}
        \centering
        \includegraphics[width = \textwidth]{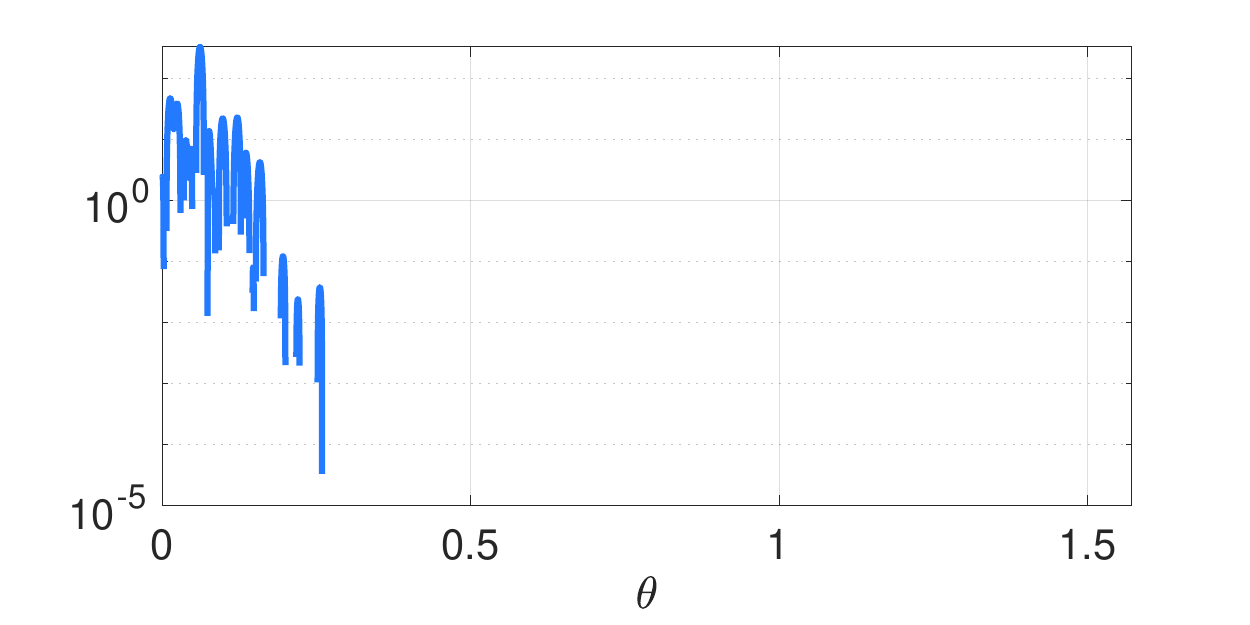}
    \end{minipage}
    \begin{minipage}[t]{0.32\textwidth}
        \centering
        \includegraphics[width = \textwidth]{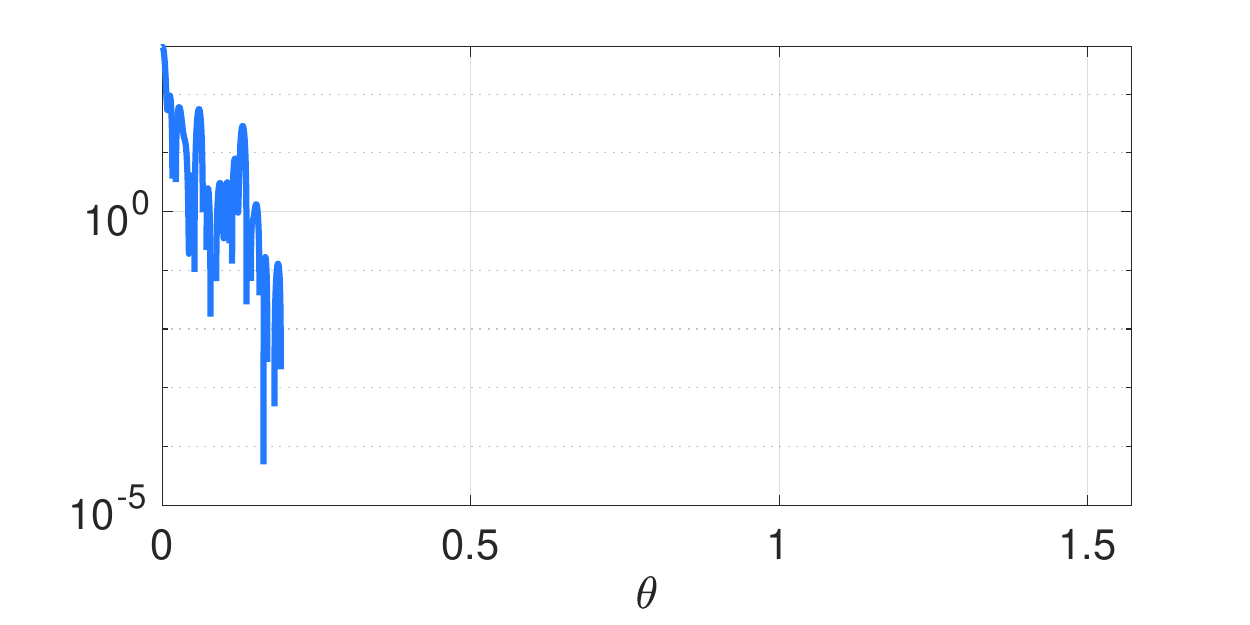}
    \end{minipage}
     \begin{minipage}[t]{0.015\textwidth}
        \rotatebox[origin=left]{90}{\ \ \ Weighted}
    \end{minipage}
    \begin{minipage}[t]{0.32\textwidth}
        \centering
        \includegraphics[width = \textwidth]{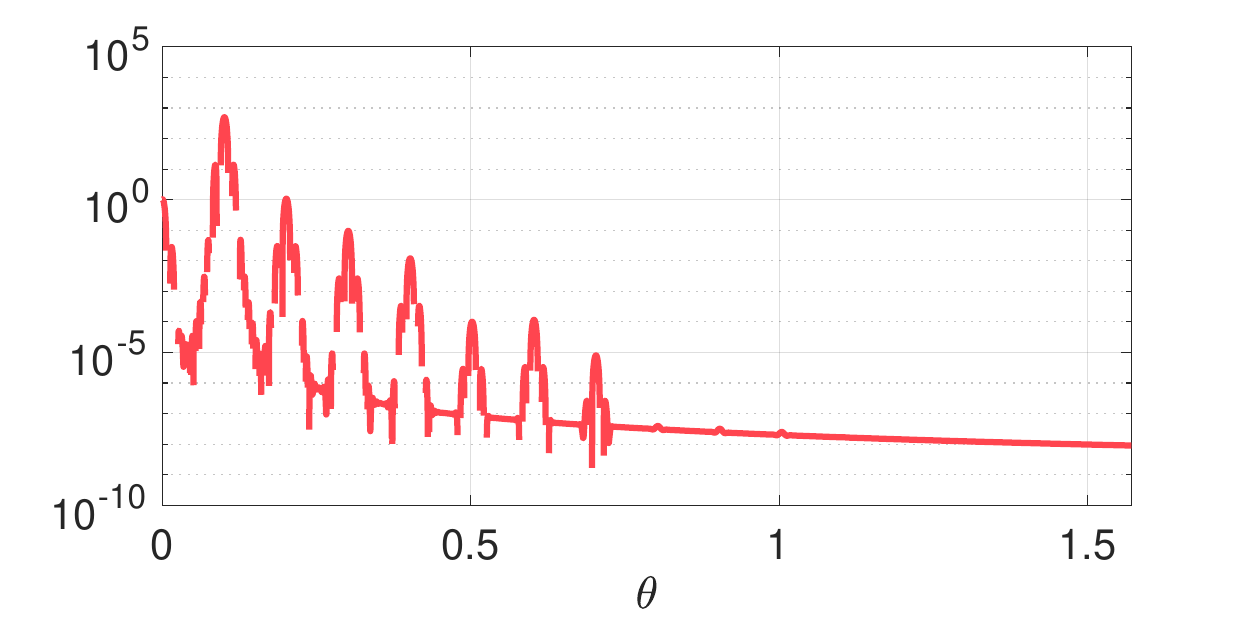}
    \end{minipage}
    \begin{minipage}[t]{0.32\textwidth}
        \centering
        \includegraphics[width = \textwidth]{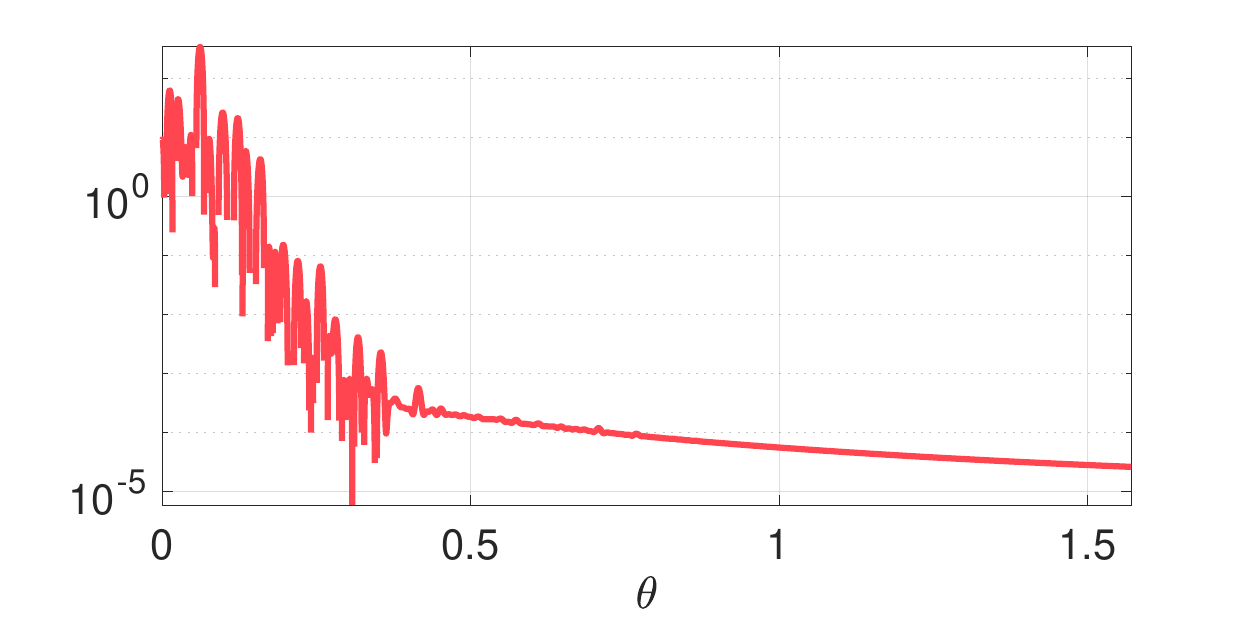}
    \end{minipage}
    \begin{minipage}[t]{0.32\textwidth}
        \centering
        \includegraphics[width = \textwidth]{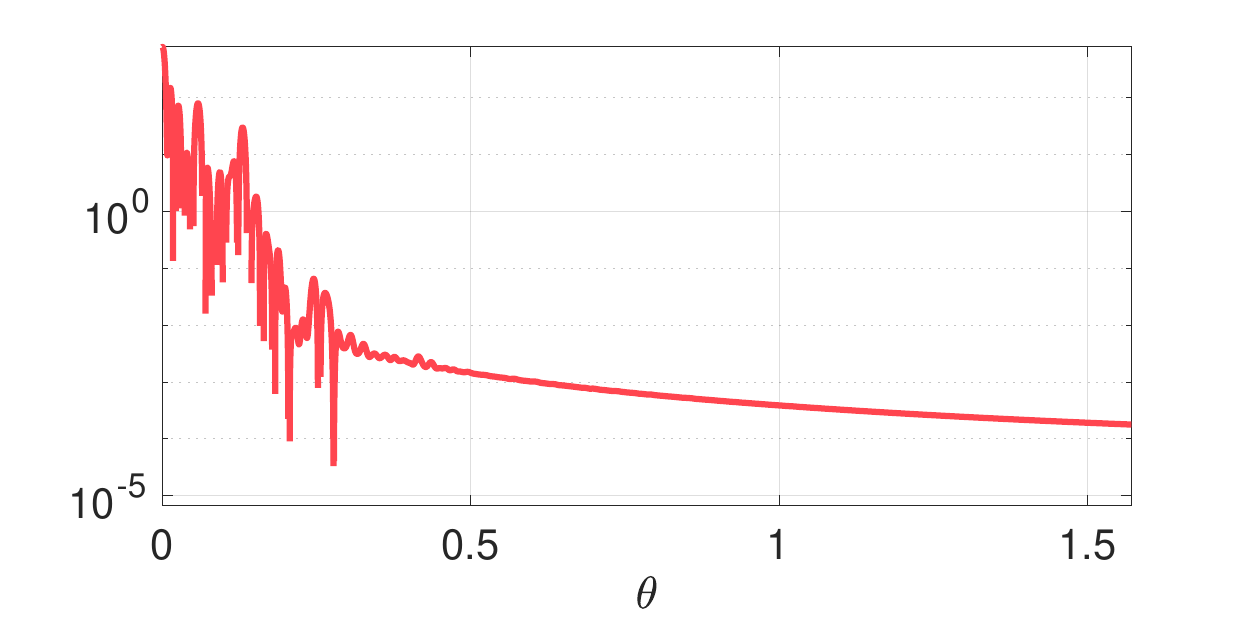}
    \end{minipage}
    \caption{Comparison between computed density approximations \eqref{spectraldensity} using $M = 1000$ unweighted (top) and weighted (bottom) autocorrelations and only $N = 1100$ data snapshots. Reynolds numbers are $13,000$ (left), $16,000$ (middle), and $19,000$ (right).} 
    \label{fig:specmeasures}
\end{figure}

We compute the first $M = 1000$ autocorrelations \eqref{autocorrelation} using both unweighted and weighted averages, from only the first $N$ snapshots of the kinetic energy. Figure~\ref{fig:auto_error} shows the relative error of these autocorrelation estimates compared with those obtained using all $N = 20,000$ snapshots.
As a further test, we use the
1000 unweighted and weighted autocorrelations computed from only 1100 snapshots (5.5\% of the total) to form the spectral density approximation \eqref{spectraldensity}.
Here we use a sharp cosine-type filter
\begin{equation}
    \varphi(x) = \frac{1}{2} + \frac{1}{2}\cos(\pi x), \quad p(x) = 35x^4 - 84x^6 - 20 x^7.  
\end{equation}
Figure~\ref{fig:specmeasures} plots the resulting densities for $\theta \in [0, \frac \pi 2]$, comparing unweighted (top row) and weighted (bottom row) constructions. In the periodic case ($\mathrm{Re}=13,000$), the weighted density cleanly identifies a set of sharp peaks at harmonics of the base shedding frequency.
These peaks are faint, merged, or even absent in the unweighted estimate. At higher Reynolds numbers, the two methods agree near $\theta \approx 0$, but differ at larger $\theta$. In particular, the unweighted density often goes negative for $\theta \gtrsim 0.3$, and so disappears on the logarithmic vertical scale, whereas the weighted density remains well resolved. This is exactly the kind of improvement we seek: clearer spectral structure from far less data.

\begin{algorithm}[t]
    \caption{Weighted measure-preserving EDMD. {\color{black}The use of QR decompositions differs from \cite{colbrook2023mpedmd} and improves numerical robustness \cite{colbrook2024multiverse}.}}\label{alg:mpEDMD}
    {\bf Input:} Sequential snapshots $\{X_n\}_{n = 1}^{N+1}$, bump function $w \in \mathcal{W}$, dictionary of observables $\vec\psi = (\psi_1\, \cdots \, \psi_L)$.
    \begin{enumerate}
        \item Compute the matrices $\Psi$ and $\Phi$ defined in \eqref{DictionaryData} with $\vec\phi = \vec\psi$ and set $W = \mathrm{diag}(w(0),w(1/N),\dots,w((N-1)/N))$.
        \item {\color{black}Compute a column-pivoted economy QR decomposition $W^{1/2}\Psi=QRP^\top$.\\
        \noindent
        ($Q$ has orthonormal columns, $R$ is upper triangular with positive diagonals, and $P$ is a permutation matrix)}
        \item {\color{black}Compute an SVD of $(PR^{-1})^{*}\Phi^*W^{1/2}Q=U_1\Sigma U_2^*$.}
        \item Compute the eigendecomposition of the unitary matrix $U_2U_1^*=\hat{V}\Lambda \hat{V}^*$.\\
        \noindent(via a Schur decomposition)
        \item {\color{black}Compute $K^{(N)}_w =PR^{-1}U_2U_1^*RP^\top$ and $V=PR^{-1}\hat{V}$.}
    \end{enumerate}
    {\bf Output:} Approximate Koopman matrix $\mathcal{K}^{(N)}_w$ with eigenvectors $V$ and eigenvalues $\Lambda$.
\end{algorithm}

\begin{figure}[t] 
    \center
    \begin{minipage}[t]{0.015\textwidth}
        \rotatebox[origin=left]{90}{\ \ \ \ \ \ \ \ \ \ \ \ \ \ \  Weighted \ \ \ \ \ \ \ \ \ \ \ \ \ \ \ \ \ Unweighted}
    \end{minipage}
    \begin{minipage}[t]{0.975\textwidth}
        \centering
        \includegraphics[width = \textwidth]{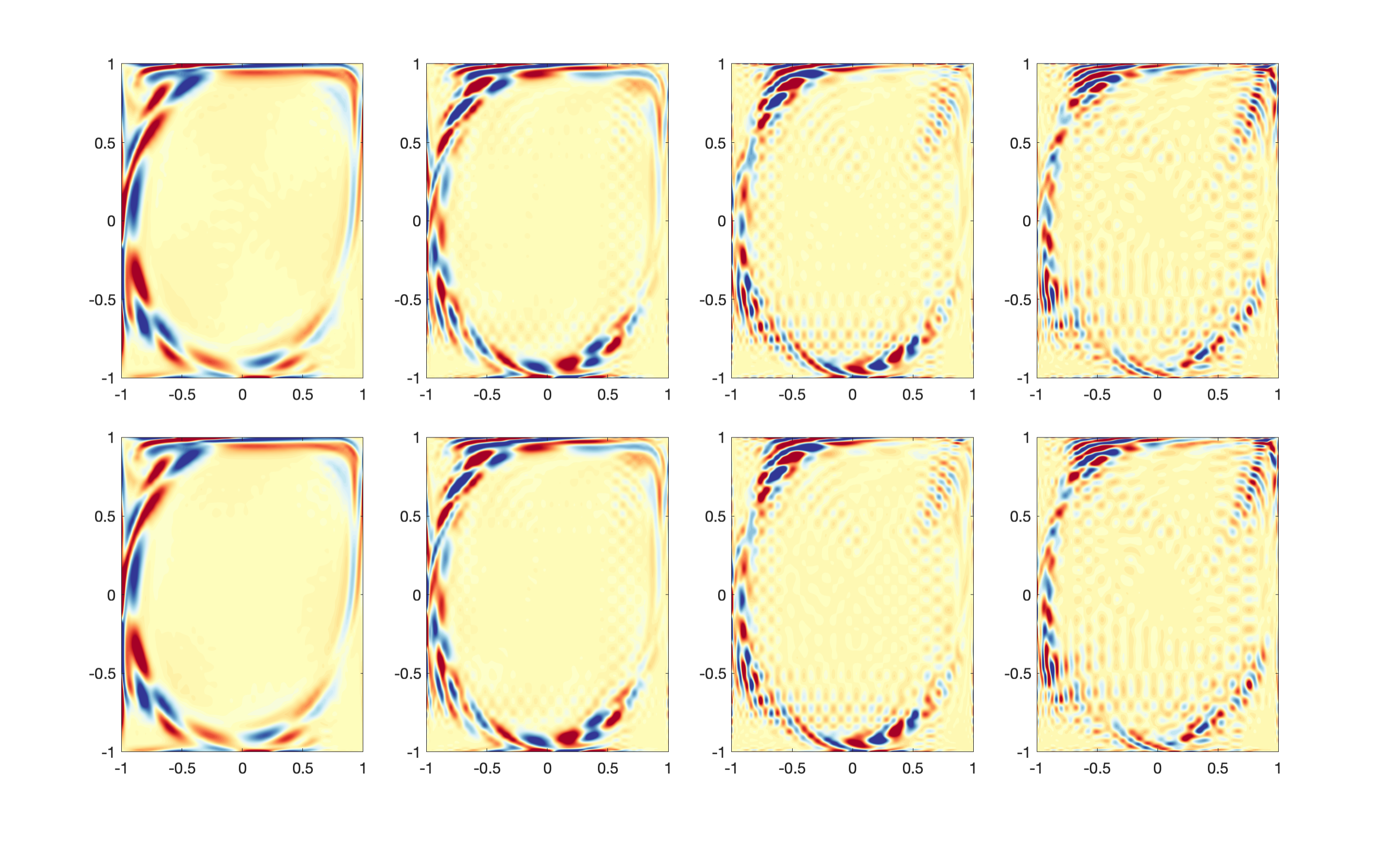}
    \end{minipage}
    \caption{Generalised eigenfunctions of lid-driven cavity flow at Reynolds number $13,000$ corresponding to $\theta = \theta_\mathrm{per},2\theta_\mathrm{per},3\theta_\mathrm{per},4\theta_\mathrm{per}$ (left to right) in the weighted density approximation. Computations use the unweighted (top row) and weighted (bottom row) rigged DMD algorithm and $N = 1100$ snapshots.}
    \label{fig:rigged_13k_1to4}
\end{figure}

\begin{figure} 
    \center
    \begin{minipage}[t]{0.015\textwidth}
        \rotatebox[origin=left]{90}{\ \ \ \ \ \ \ \ \ \ \ \ \ \ \  Weighted \ \ \ \ \ \ \ \ \ \ \ \ \ \ \ \ \ Unweighted}
    \end{minipage}
    \begin{minipage}[t]{0.975\textwidth}
        \centering
        \includegraphics[width = \textwidth]{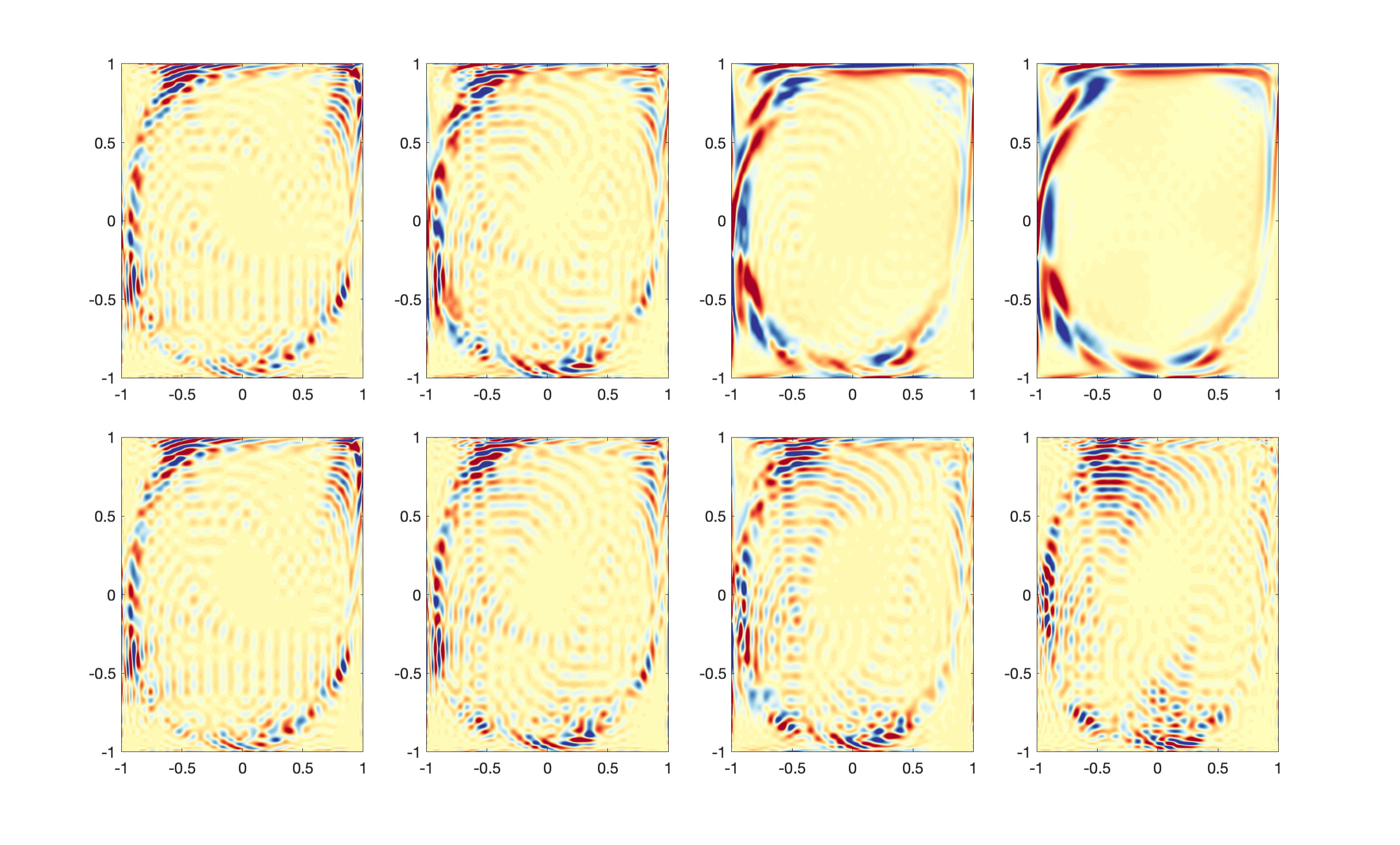}
    \end{minipage}
    \caption{Generalized eigenfunctions of lid-driven cavity flow at Reynolds number $13,000$ corresponding to $\theta = 5\theta_\mathrm{per},6\theta_\mathrm{per},7\theta_\mathrm{per},8\theta_\mathrm{per}$ (left to right) in the weighted density approximation. Computations use the unweighted (top row) and weighted (bottom row) rigged DMD algorithm and $N = 1100$ snapshots.}
    \label{fig:rigged_13k_5to8}
\end{figure}

\begin{figure} 
    \center
    \begin{minipage}[t]{0.015\textwidth}
        \rotatebox[origin=left]{90}{\ \ \ \ \ \ \ \ \ \ \ \ \ \ \  Weighted \ \ \ \ \ \ \ \ \ \ \ \ \ \ \ \ \ Unweighted}
    \end{minipage}
    \begin{minipage}[t]{0.975\textwidth}
        \centering
        \includegraphics[width = \textwidth]{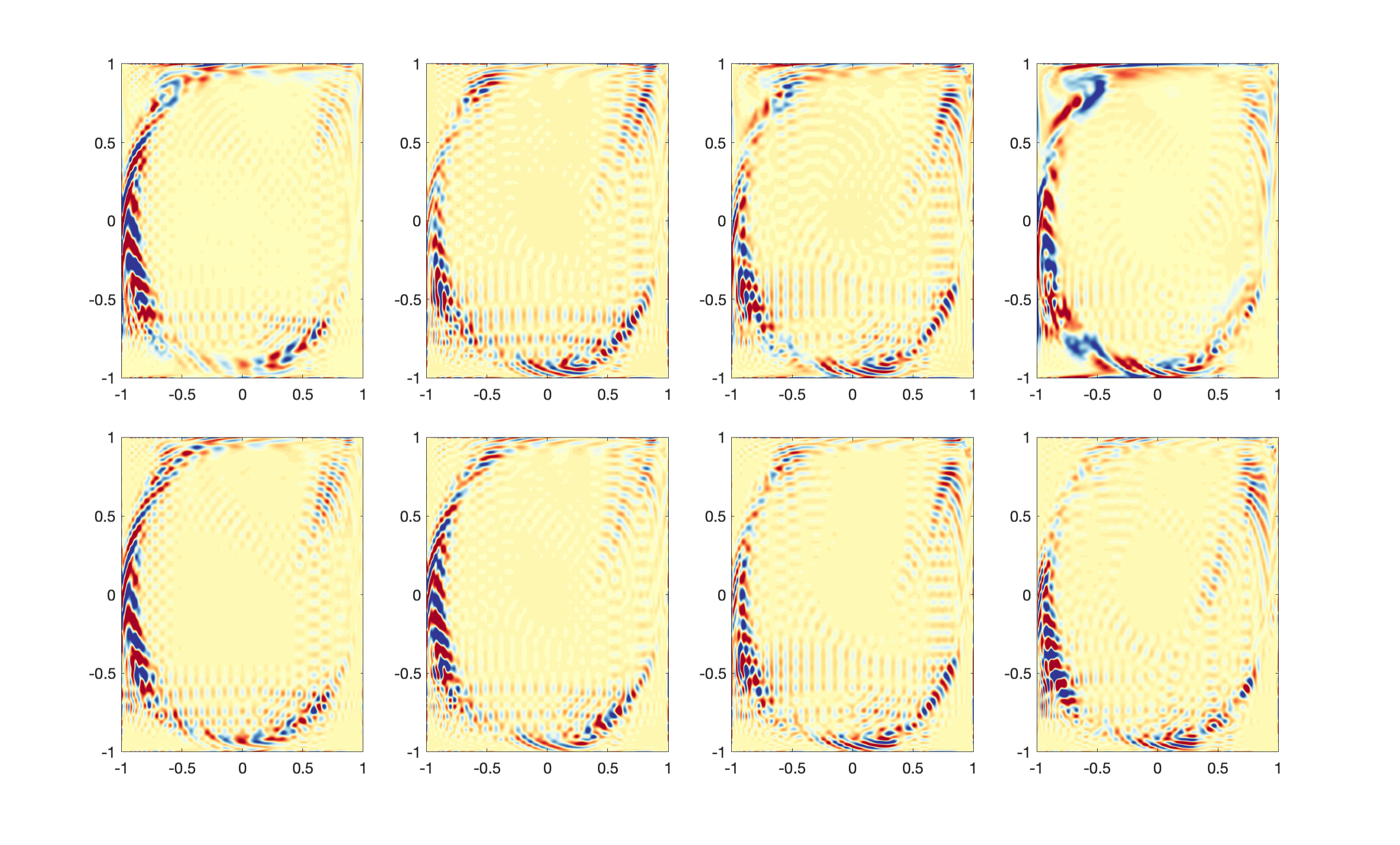}
    \end{minipage}
    \caption{The same as Figure~\ref{fig:rigged_13k_5to8} but with Reynolds number $16,000$.}
    \label{fig:rigged_16k_5to8}
\end{figure}

We can also use these spectral approximations in rigged DMD \cite{colbrook2025rigged} to identify generalised Koopman eigenfunctions.
A key ingredient is measure-preserving EDMD \cite{colbrook2023mpedmd}, which can likewise incorporate weights; see Algorithm~\ref{alg:mpEDMD}.
We apply this to the cavity flow data, using only $N = 1100$ snapshots.
We extract approximate generalised eigenfunctions at angles $\theta$ corresponding to multiples of the base frequency in the periodic flow at $\mathrm{Re}=13,000$, denoted $\theta_\mathrm{per} \approx 0.1$.
Figure~\ref{fig:rigged_13k_1to4} shows the resulting modes for $\theta \in \{\theta_\mathrm{per},2\theta_\mathrm{per},3\theta_\mathrm{per},4\theta_\mathrm{per}\}$.
Here, weighted and unweighted methods give broadly similar structures.

For higher harmonics, the distinction becomes clearer.
Figure~\ref{fig:rigged_13k_5to8} shows the modes associated with $\theta \in \{5\theta_\mathrm{per},6\theta_\mathrm{per},7\theta_\mathrm{per},8\theta_\mathrm{per}\}$.
The unweighted method often returns modes with little discernible spatial structure, whereas the weighted method still recovers meaningful coherent patterns. A similar effect is visible at 
$\mathrm{Re}=16,000$ in Figure~\ref{fig:rigged_16k_5to8}: with limited data, the weighted procedure continues to resolve generalised eigenfunctions at higher frequencies.

\subsection{Example Method V: Diffusion forecasts}

Our final demonstration of weighted averages concerns nonparametric forecasting.
Diffusion forecasting methods \cite{giannakis2019data,berry2015nonparametric,thiede2019galerkin,giannakis2020extraction,berry2015nonparametric2} build predictive models directly from data without first fitting a parametric surrogate model.
These methods all rely on diffusion maps \cite{coifman2005geometric,coifman2006diffusion} to construct an effective basis of observables from sampled data.

Following \cite{berry2015nonparametric}, we assume that the equally spaced samples $\{X_n\}_{n = 1}^{N+1}$ come from an underlying stochastic differential equation
\begin{equation}\label{SDE}
    \mathrm{d}x = a(x)\,\mathrm{d}t + b(x)\,\mathrm{d}B_t,
\end{equation}
on a manifold $\mathcal{X} \subset \mathbb{R}^d$, where $B_t$ is standard Brownian motion, $a$ is a drift vector field, and $b$ is a diffusion tensor. The aim is forecasting.
Given an initial probability density $p_0(x)$, we want to estimate the evolved density $p(x,t) = \mathrm{e}^{t\mathcal{L}^*}p_0(x)$ for $t > 0$, where $\mathcal{L}^*$ is the Fokker–Planck operator of \eqref{SDE},
without explicitly estimating $a$, $b$, or $\mathcal{L}^*$.

We assume that the data $\{X_n\}_{n = 1}^{N+1}$ are sampled from an invariant (equilibrium) density $p_\mathrm{eq}(x)$ of \eqref{SDE}, i.e., $p_\mathrm{eq}(x) = \mathrm{e}^{t\mathcal{L}^*}p_\mathrm{eq}(x)$ for all $t > 0$. Using these data, one applies a variable-bandwidth diffusion maps construction \cite{berry2016variable} to learn eigenfunctions $\phi_j$ that are orthonormal in $L^2(p_\mathrm{eq}^{-1})$. The diffusion forecast \cite{berry2015nonparametric} for the density one sampling interval $\tau$ in the future is then
\begin{equation}\label{DiffusionForecast}
    p(x,\tau) \approx \sum_{i = 1}^M p_\mathrm{eq}(x)\phi_i(x)\sum_{j = 1}^M A_{ij}c_j
\end{equation}
where $\tau > 0$ is the data sampling rate, $M \geq 1$ is the number of computed eigenfunctions $\{\varphi_m\}_{m = 1}^M$, $c_j = \langle p_0,\phi_j\rangle$, and 
\begin{equation}\label{DiffusionAverage}
    A_{ij} = \frac{1}{N}\sum_{n = 1}^N \phi_j(X_n)\phi_i(X_{n+1}).
\end{equation}
The forward-time approximation of $p(x,k\tau)$ with $k \geq 1$ is achieved by replacing $A_{ij}$ in \eqref{DiffusionForecast} with the $(i,j)$th component of the matrix power $A^k$, where $A = [A_{ij}]_{i,j = 1}^M$. Once again, we recognise an ergodic average in \eqref{DiffusionAverage}. Once again, we may replace it by a weighted average using a taper $w \in \mathcal{W}$, yielding a weighted diffusion forecast. 

\paragraph{Application to El-Ni\~{n}o dataset.}
We now compare standard and weighted diffusion forecasts on the Niño-3.4 index, which measures monthly sea-surface temperature anomalies in a region of the equatorial Pacific commonly used to characterise El Niño/La Niña variability.
The goal here is not to compete with state-of-the-art climate models, but simply to compare weighted versus unweighted estimation of $A_{ij}$ in \eqref{DiffusionAverage}.

We follow the setup of \cite{berry2015nonparametric}.
We construct a 6-lag time-delay embedding of the Niño-3.4 index and compute the diffusion map eigenfunctions.
We train the model on monthly data from January 1920 through December 1999, and we evaluate forecasts on January 2000–December 2013.
We use $M = 14$ eigenfunctions in \eqref{DiffusionForecast}.

\begin{figure}[t] 
    \center
    \begin{minipage}[t]{0.48\textwidth}
        \centering
        Unweighted Diffusion Forecast
    \end{minipage}
    \begin{minipage}[t]{0.48\textwidth}
        \centering
        Weighted Diffusion Forecast
    \end{minipage}
    \begin{minipage}[t]{0.48\textwidth}
        \centering
        \includegraphics[width = \textwidth]{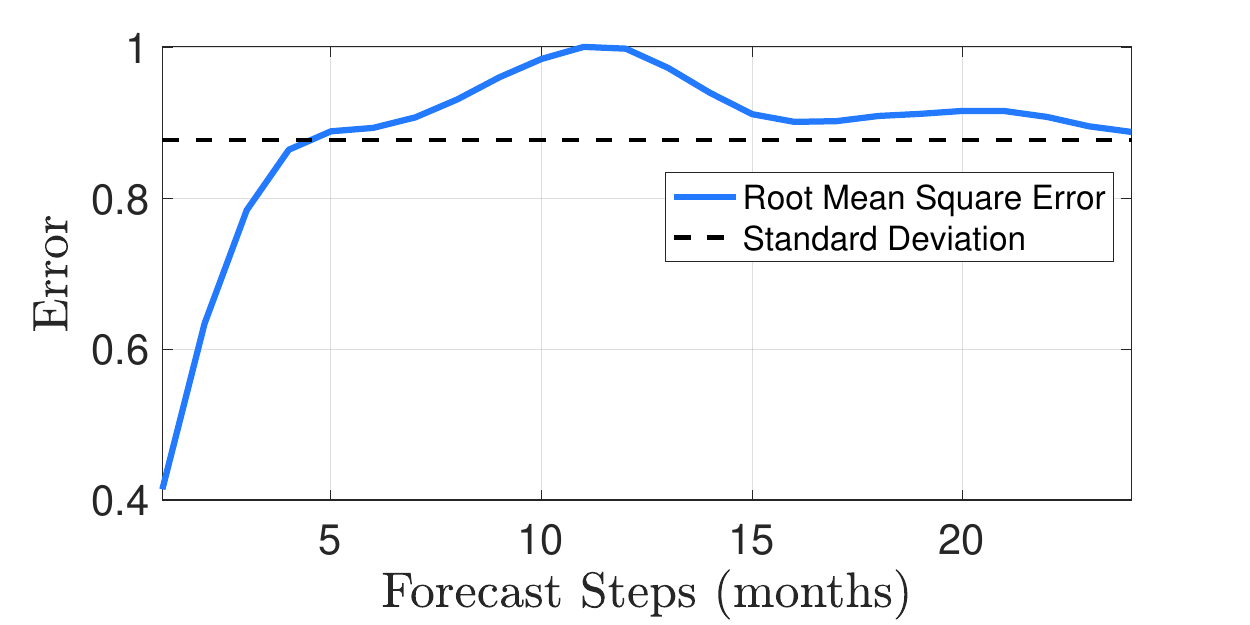}
    \end{minipage}
    \begin{minipage}[t]{0.48\textwidth}
        \centering
        \includegraphics[width = \textwidth]{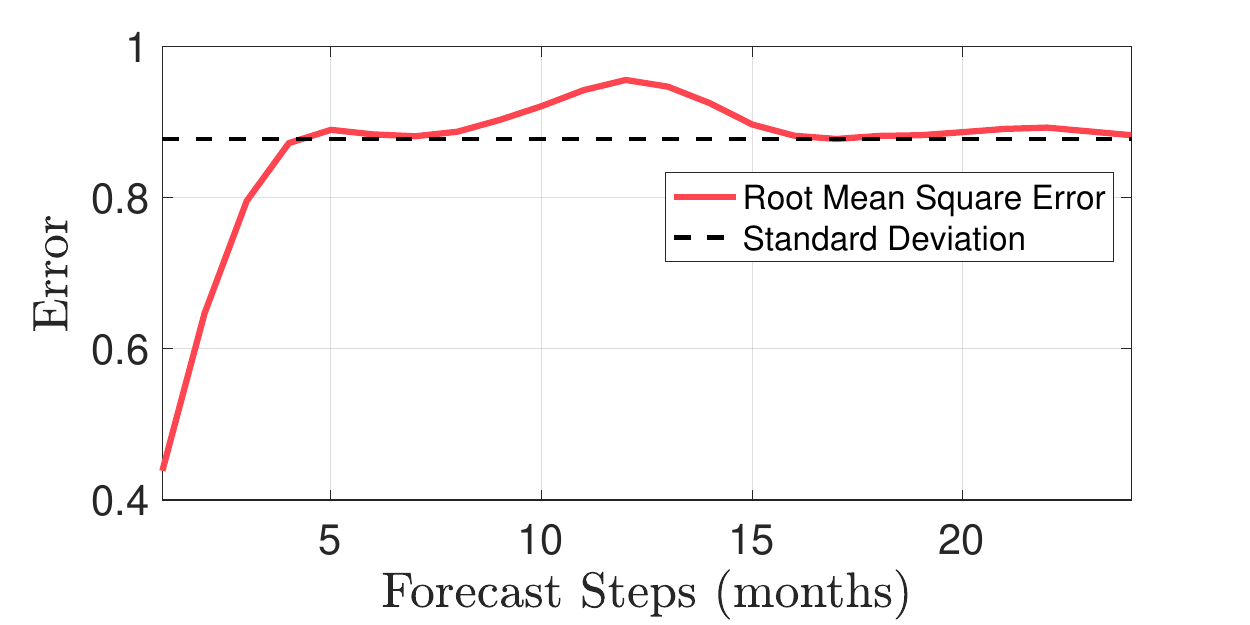}
    \end{minipage}
    \begin{minipage}[t]{0.48\textwidth}
        \centering
        \includegraphics[width = \textwidth]{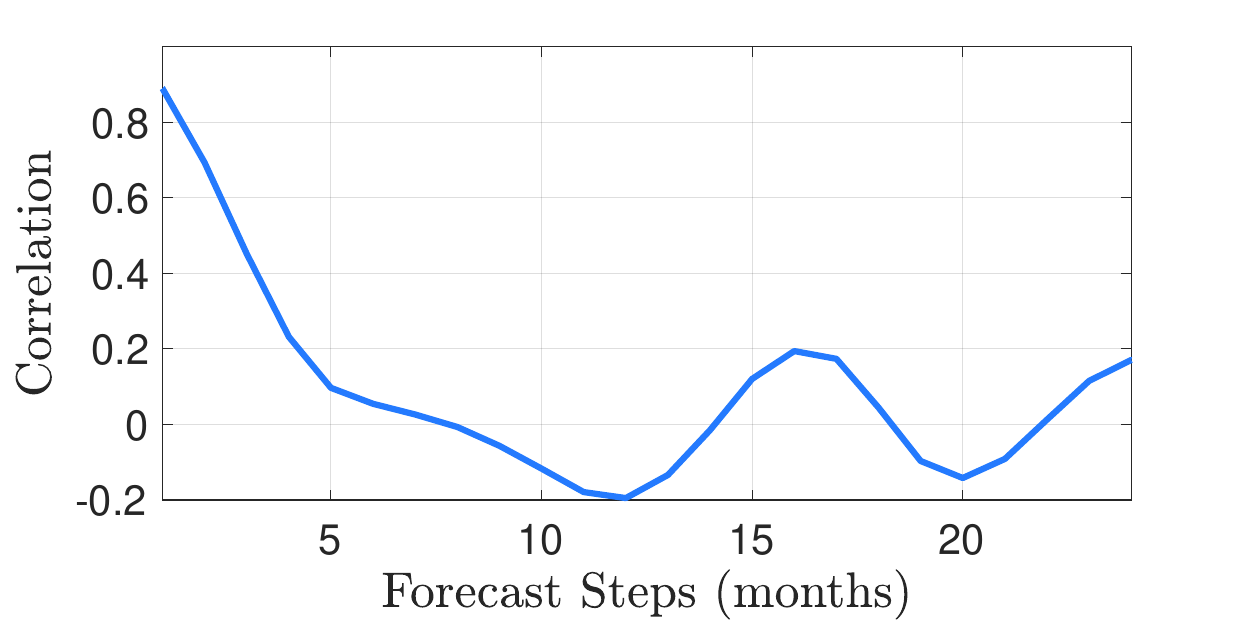}
    \end{minipage}
    \begin{minipage}[t]{0.48\textwidth}
        \centering
        \includegraphics[width = \textwidth]{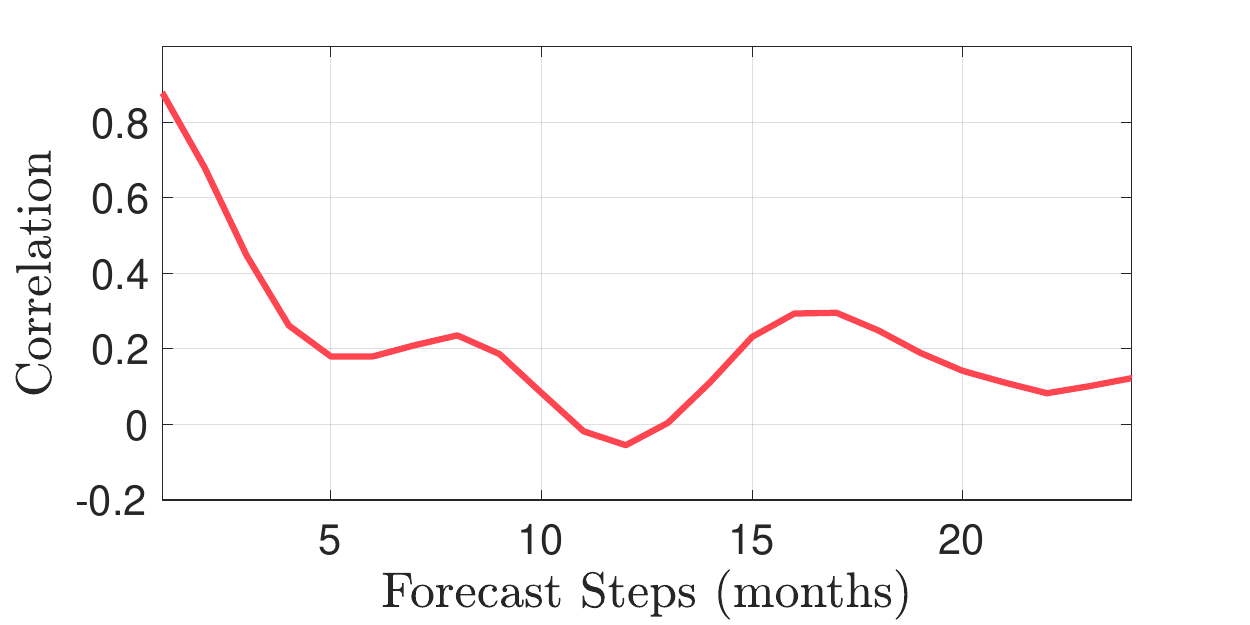}
    \end{minipage}
    \begin{minipage}[t]{0.48\textwidth}
        \centering
        \includegraphics[width = \textwidth]{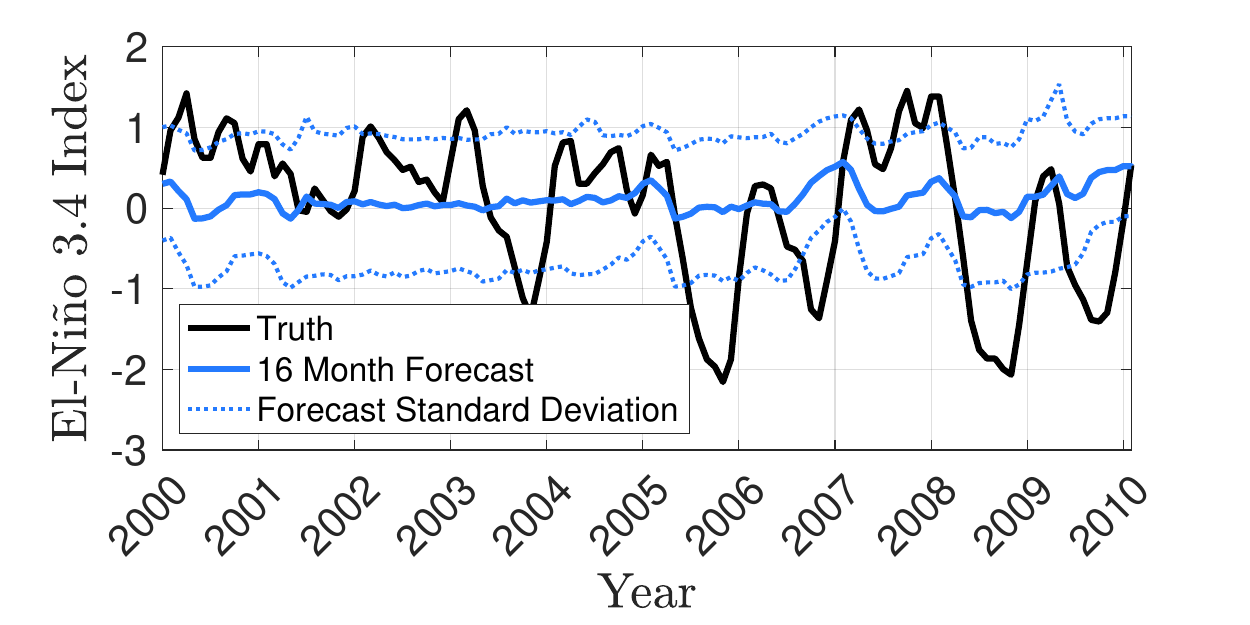}
    \end{minipage}
    \begin{minipage}[t]{0.48\textwidth}
        \centering
        \includegraphics[width = \textwidth]{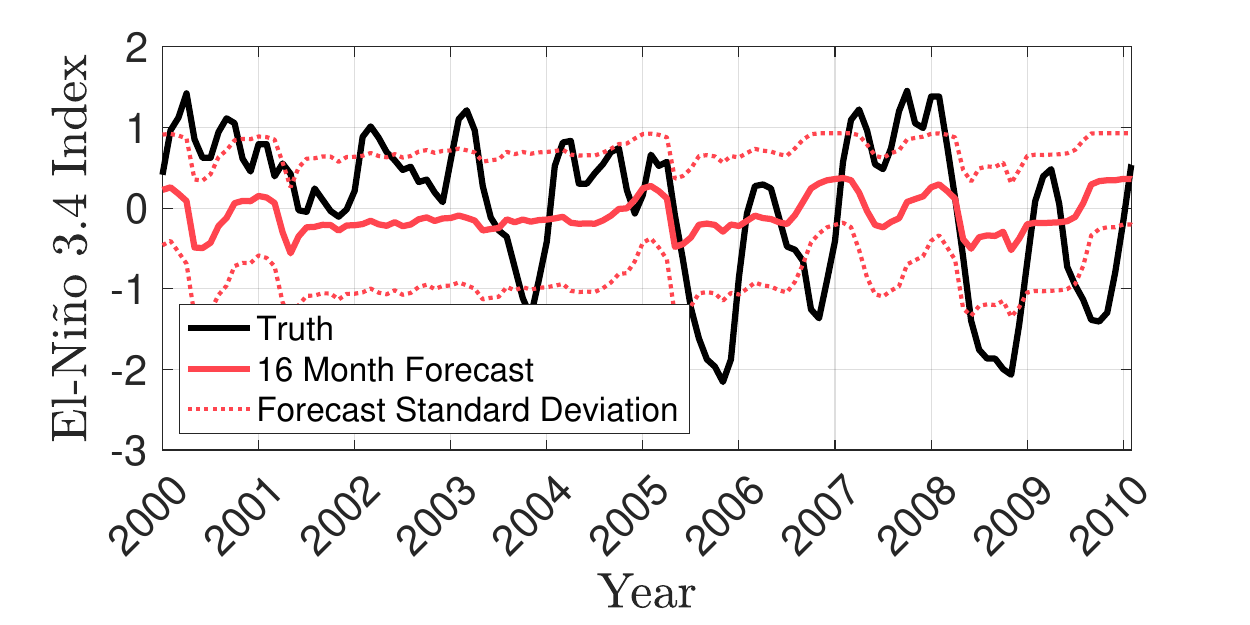}
    \end{minipage}
    \caption{A comparison of unweighted (left) and weighted (right) diffusion forecasts for the El-Ni\~no 3.4 index using $M = 14$ diffusion map eigenfunctions.} 
    \label{fig:diffusionforecast}
\end{figure}

Figure~\ref{fig:diffusionforecast} compares unweighted (left column) and weighted (right column) diffusion forecasts.
The first row shows the root mean square error (RMSE) of the forecast as a function of lead time.
Both forecasts degrade to the climatological error (the standard deviation of the validation data) after about 5 months of lead time.
Interestingly, the error improves again near a lead time of about 16 months, and the weighted forecast stays at or slightly below the climatological RMSE there.
The second row shows correlation skill versus lead time.
Both correlations fall off and then recover around the 16-month lead time, with the weighted forecast peaking about 0.1 higher.
The unweighted forecast correlation becomes negative for lead times of roughly 7–14 months;
the weighted forecast correlation only dips negative at 11 and 12 months.
The bottom row plots the actual 16-month lead forecast time series for both methods.

Similar results are returned for $10 \leq M \leq 30$ eigenfunctions. For $M < 10$, the basis is too small: the expansion in \eqref{DiffusionForecast} is not expressive enough to capture the observed variability.
For $10 \leq M \leq 30$, the overall behaviour is similar to the $M=14$  case.

For much larger $M$, the picture worsens.
Using a large number of eigenfunctions substantially degrades the weighted forecast compared with the unweighted one.
We do not have a definitive explanation for this, but a plausible cause is that higher-index diffusion map eigenfunctions are poorly resolved from the available data ($N = 948$ monthly samples in the training set).
Weighted Birkhoff averages only accelerate convergence cleanly for smooth observables (see Section~\ref{sec:WBAs}).
If the high-order eigenfunctions are inaccurate or noisy, their weighted averages may amplify that error rather than average it out.

Despite that caveat, the weighted diffusion forecast with $M=14$ achieves consistently lower RMSE and higher correlation, and does so across a wide range of lead times, while using the same data and the same diffusion-map basis as the unweighted method.

\section{Discussion}

In this work, we have shown how weighted Birkhoff averages can be built into a range of data-driven methods to extract more from limited data. We focused on five widely used tools for analysing dynamical systems: Dynamic Mode Decomposition, Koopman operator approximations, model identification, spectral measures, and diffusion forecasting. In each case, weighting can be introduced with only minor changes to standard practice.

The weighted variants were tested on canonical examples, such as flow past a cylinder and lid-driven cavity flow, as well as on real data from the Niño-3.4 index and on reduced-order models with noisy numerical derivatives.
Our aim has been expository: the same idea can be applied to many other methods that depend on long-time averages.

Although we have presented examples rather than theorems, convergence guarantees for the weighted versions follow directly from established results on weighted Birkhoff averages.
More interesting, perhaps, is what happens beyond the theory.
Even in settings where one would not expect accelerated convergence---chaotic dynamics, stochastic forcing, measurement noise---the weighted methods often match or outperform their unweighted counterparts.
And even when they do not, they typically produce smoother and more stable estimates as the data length varies (see Figures \ref{fig:Logistic}, \ref{fig:wSINDy}, and \ref{fig:auto_error}), improving reliability in practice.

Weighted averages are easy to implement, inexpensive to compute, and often better.
They offer a practical route to more accurate operator learning, forecasting, and system identification, not by adding new machinery, but by making better use of the data one already has.

\section*{Acknowledgments}

The authors thank Evelyn Sander for a helpful discussion explaining the utility of weighted Birkhoff averages that initiated this project and Tyrus Berry for providing code and insight for the diffusion forecast example. M.B.S.E.T. was supported by a Concordia Undergraduate Summer Research Award. J.J.B. was partially supported by the Natural Sciences and Engineering Research Council of Canada (NSERC) through grant RGPIN-2023-04244 and the Fonds de Recherche du Qu\'ebec Nature et Technologies (FRQNT) through the grant 340894. M.J.C. was supported by Isaac Newton Trust Grant No. LEAG/929.

\section*{Data Availability Statement}

All code and data required to reproduce the results in this paper can be found in the repository \href{https://github.com/jbramburger/weighted_methods}{https://github.com/jbramburger/weighted\_methods}.

\footnotesize
\renewcommand{\baselinestretch}{0.95}
\bibliographystyle{abbrv}
\bibliography{references.bib}

\end{document}